\documentclass[11pt,a4paper]{article}

\usepackage[T1]{fontenc}
\usepackage{lmodern}
\usepackage[margin=28mm]{geometry}
\usepackage{amsmath,amssymb,amsthm,mathtools}
\usepackage{microtype}
\usepackage{needspace}
\usepackage{tikz}
\usetikzlibrary{arrows.meta}
\usepackage[font=small,labelfont=bf]{caption}

\definecolor{ink}{HTML}{24384B}
\definecolor{sea}{HTML}{007F83}
\definecolor{indigo}{HTML}{4D609C}
\definecolor{ochre}{HTML}{B26820}
\usepackage[colorlinks=true,linkcolor=ink,citecolor=indigo,urlcolor=indigo]{hyperref}
\hypersetup{pdftitle={Improved Lower Bound on the Hadwiger-Debrunner Numbers}}

\newif\ifshowchanges
\showchangesfalse
\definecolor{revisionred}{HTML}{A51E2D}
\definecolor{revisionblue}{HTML}{164CA4}
\makeatletter
\newcommand{\RevisionInk}[1]{%
  \color{#1}%
  \def\@linkcolor{#1}%
  \def\@citecolor{#1}%
  \def\@urlcolor{#1}%
}
\AtBeginDocument{\let\RevisionMakeFnText\@makefntext}
\newcommand{\RevisionColor}[1]{%
  \RevisionInk{#1}%
  \long\def\@makefntext##1{%
    \RevisionInk{#1}%
    \RevisionMakeFnText{##1}%
  }%
}
\makeatother
\DeclareRobustCommand{\added}[1]{%
  \ifshowchanges{\RevisionColor{revisionblue}#1}\else#1\fi}
\DeclareRobustCommand{\deleted}[1]{%
  \ifshowchanges{\RevisionColor{revisionred}#1}\fi}
\newenvironment{additions}
  {\ifshowchanges\RevisionColor{revisionblue}\fi}
  {}
\newtheorem{theorem}{Theorem}[section]
\newtheorem{proposition}[theorem]{Proposition}
\newtheorem{lemma}[theorem]{Lemma}

\newcommand{\R}{\mathbb R}
\newcommand{\Prob}{\mathbb P}
\newcommand{\one}{\mathbf 1}
\newcommand{\HD}{\operatorname{HD}}
\newcommand{\depth}{\operatorname{depth}}
\newcommand{\area}{\operatorname{area}}
\newcommand{\diff}{\,\mathrm d}

\tikzset{
  axis/.style={draw=ink!75,semithick,-{Latex[length=2mm]}},
  strip a/.style={draw=sea,semithick},
  strip b/.style={draw=indigo,semithick,densely dashed},
  strip c/.style={draw=ochre,semithick,dash pattern=on 3pt off 1.5pt on .6pt off 1.5pt}
}

\title{Improved Lower Bound on the Hadwiger-Debrunner Numbers}

\author{
Chaya Keller\thanks{Institute for the Theory of Computing,
Faculty of Computer and Information Science,
Ben-Gurion University of the Negev,
Be'er Sheva 84105, Israel and
School of Computer Science, Ariel University, Israel. \texttt{chayak@bgu.ac.il}
}}

\date{}

\begin{document}
\maketitle

\begin{abstract}
The celebrated $(p,q)$-theorem of Alon and Kleitman asserts that for any $d \in \mathbb{N}$ and any $p \geq q \geq d+1$, every family of convex sets in $\mathbb{R}^d$ in which among any $p$ sets, some $q$ intersect, can be pierced by $c(p,q,d)$ points. The minimum numbers $c(p,q,d)$ are called the Hadwiger-Debrunner numbers and denoted $\mathrm{HD}_d(p,q)$ and determining them is a longstanding open problem. We prove that for any $d \geq 2$, $\mathrm{HD}_d(p,d+1) \geq p^{d-o(1)}$, and more generally, $\mathrm{HD}_d(p,q) \geq p^{\frac{q-1}{q-d}-o(1)}$. Our bound is obtained by a family of parallelotopes, and is tight up to polylogarithmic factors among all families that have a bounded VC dimension. The proof uses the hypergraph container method.
\end{abstract}

\section{Introduction}

\subsection{Background}

\textbf{The $(p,q)$ theorem and the Hadwiger-Debrunner numbers.} For $p,q \in \mathbb{N}$, a family $\mathcal{F}$ of sets is said to satisfy the $(p,q)$ property if among any $p$ sets in $\mathcal{F}$, some $q$ intersect. In this notation, the classical Helly's theorem says that a finite family of at least $d+1$ convex sets in $\R^d$ that satisfies the $(d+1,d+1)$ property can be pierced by a single point. In 1957, Hadwiger and Debrunner~\cite{HD57} conjectured that for any $p \geq q \geq d+1$, there exists $c(p,q,d)$ such that any finite family of at least $p$ convex sets in $\R^d$ that satisfies the $(p,q)$ property can be pierced by $c(p,q,d)$ points. In their classical $(p,q)$-theorem from 1992, Alon and Kleitman~\cite{AK92} proved that $c(p,q,d)$ is indeed finite for any $p \geq q\geq d+1$. The minimum piercing numbers $c(p,q,d)$ are called the \emph{Hadwiger-Debrunner numbers} (or in short, \emph{HD numbers}), and are denoted by $\HD_d(p,q)$. 

\medskip \noindent \textbf{Quantitative bounds on the HD numbers.} The upper bound on $\HD_d(p,q)$ in the proof of Alon and Kleitman is very large. They discussed possible ways to improve it and concluded their paper stating the problem of estimating the numbers $\HD_d(p,q)$ more precisely. This question has become one of the best-known quantitative Helly-type problems (see, e.g.,~\cite[Problem~2.12]{BK22}) and has been studied extensively in the last 30 years, both in general (e.g.,~\cite{BMN11,KellerS21,KST18,RN26}) and in special cases (e.g.,~\cite{KGT01,MG22}); yet, it remains wide open. The most well-studied case is $q=d+1$, for which the $(p,q)$ property assumption is the weakest (e.g.,~\cite{KGT01,MG22,RN26}).

The best known general upper bound, obtained by Keller, Smorodinsky and Tardos~\cite{KST18}, is $\HD_d(p,q) \leq \tilde{O}(p^{\frac{q-1}{q-d}\cdot \rho_d})$, where $\rho_d \approx d-0.5$ is the exponent in the \emph{weak $\varepsilon$-net theorem}~\cite{Rubin22,Rubin26} (see below). For $q \geq \log p$, stronger bounds were obtained in~\cite{KST18}.  
Lower bounds were obtained in several works that can be divided into two lines of work which we describe below.

\medskip \noindent \textbf{Lower bounds on the HD numbers via weak $\epsilon$-nets.}  The first line derives a lower bound on $\HD_d(p,q)$ from a lower bound in the \emph{weak $\varepsilon$-net theorem}. A weak $\varepsilon$-net for a set $P \subset \R^d$ is a set $N \subset \R^d$ such that every convex set that contains at least $\varepsilon |P|$ points of $P$, contains a point of $N$. The weak $\varepsilon$-net theorem asserts that any $P \subset \R^d$ admits a weak $\varepsilon$-net of size at most $f_d(\varepsilon)$. As was noted in~\cite{KellerS21}, for any $P \subset \R^d$ and any $p,q \in \mathbb{N}$, the family $\mathcal{F}$ of convex sets in $\R^d$ that contain at least $\varepsilon |P|$ points of $P$, for $\varepsilon=\frac{q}{p}$, satisfies the $(p,q)$ property. Hence, any construction showing that $f_d(\varepsilon)>\ell$ for $\varepsilon=q/p$, implies the lower bound $\HD_d(p,q)>\ell$. The first lower bound on $\HD_d(p,q)$ was of this form: $\HD_d(p,q) \geq \Omega(\frac{p}{q}\log^{d-1}(\frac{p}{q}))$, implied by the lower bound $f_d(\varepsilon) \geq \Omega(\frac{1}{\epsilon}\log^{d-1}(\frac{1}{\varepsilon}))$ obtained by Bukh, Matou\v{s}ek and Nivasch~\cite{BMN11}. Very recently, Bhore, Browne, and Zheng~\cite{BBZ26} obtained a vastly improved lower bound on $f_d(\varepsilon)$, which implies the bound $\HD_d(p,q)=\Omega_d((\frac{p}{q})^{\gamma_d})$, where $\gamma_d = \frac{1}{2} \lfloor \frac{\sqrt{8d+25}-1}{2} \rfloor$, which is at least $3/2$ for $d \geq 3$ and is $\sqrt{d/2}+O(1)$ as $d \to \infty$. The problem of determining $f_d(\varepsilon)$ is wide open. In any case, the best lower bound on $\HD_d(p,q)$ that can be obtained via this route is $p^{\rho_d}$, where $\rho_d$ is the best upper bound on the exponent in the weak $\varepsilon$-net theorem. As was mentioned above, $\rho_d \approx d-0.5$ was obtained by Rubin~\cite{Rubin22,Rubin26}. 

\medskip \noindent \textbf{Lower bounds on the HD numbers via families of hyperplanes.}
The second line of work obtains lower bounds on $\HD_d(p,q)$ by constructing specific families of hyperplanes in $\mathbb{R}^d$ and using the hypergraph container lemma~\cite{BMS15,ST15} to obtain the bound. The first bound of this class was obtained by Keller and Smorodinsky~\cite{KellerS21} who used a general position construction of Balogh and Solymosi~\cite{BS18} to show that $\HD_2(p,3)=\Omega_\varepsilon(p^{6/5-\varepsilon})$, and more generally, $\HD_2(p,q)=\Omega_\varepsilon(p^{1+\frac{1-\varepsilon}{4q-7}})$ for every $\varepsilon>0$.  Roche-Newton~\cite{RN26} improved the bound on $\HD_2(p,3)$ to $\Omega_\varepsilon(p^{5/4-\varepsilon})$. Keller and Smorodinsky~\cite{KS26} obtained the bound $\HD_3(p,4)\ge p^{2-o(1)}$, using a high-dimensional general position construction by Suk and Zeng~\cite{SZ26}. Very recently, Chen, Ge, and Shu~\cite{CGS26} used a randomly sparsified family of hyperplanes to show that $\HD_d(p,q)\ge p^{1+\frac{d-1}{q-1}-o(1)}$, and in particular,
$\HD_2(p,3)\ge p^{3/2-o(1)}$. In addition, they used a transfer argument from upper bounds on the sizes of sufficiently general-positioned subsets of high-dimensional integer grids to lower bounds for $\HD_d(p,q)$ to prove the bound $\HD_d(p,d+1)\ge p^{\lfloor \frac{d+1}{2} \rfloor-o(1)}$, using another family of hyperplanes.

To summarize, the best currently known lower bounds on $\HD_d(p,d+1)$ are $\HD_2(p,3) \ge p^{3/2-o(1)}$~\cite{CGS26}, $\HD_3(p,4) \ge p^{2-o(1)}$~\cite{CGS26,KS26} and $\HD_d(p,d+1) \ge p^{\lfloor \frac{d+1}{2} \rfloor-o(1)}$ for $d \geq 4$~\cite{CGS26}. For larger values of $q$, the best known lower bound is $\HD_d(p,q)=\Omega_d((\frac{p}{q})^{\gamma_d})$~\cite{BBZ26}. 

\subsection{Our results}

In this paper we prove a new lower bound on the Hadwiger-Debrunner numbers.
\begin{theorem}\label{thm:main}
For every fixed integer $d\ge 2$, there is a constant $c_d>0$ such that, for every sufficiently large integer $p$,
\[
  \HD_d(p,d+1)\ge
  c_d\,\frac{p^d\log\log p}{(\log p)^{\max\{3d+1,d^2\}}}.
\]
In particular, $\HD_d(p,d+1)\ge p^{d-o(1)}$.  Furthermore, for any $q \geq d+2$,
\[
  \HD_d(p,q)\ge
  c_{d,q}\,\frac{p^{\frac{q-1}{q-d}}\log\log p}{(\log p)^{3d+1}}.
\]
In particular, $\HD_d(p,q)\ge p^{\frac{q-1}{q-d}-o(1)}$. The convex sets in the construction are compact parallelotopes.
\end{theorem}
\noindent \textbf{Discussion.} Theorem~\ref{thm:main} significantly improves the best previous lower bounds on $\HD_d(p,d+1)$, for any $p,d$. Furthermore, it is sharp up to polylogarithmic factors within the class of families that have a \emph{bounded VC dimension}~\cite{VC71}. Indeed, by the \emph{$\varepsilon$-net theorem}~\cite{HW87}, every such family admits an $\varepsilon$-net (and in particular, a weak $\varepsilon$-net) of size $O(\frac{1}{\varepsilon}\log(\frac{1}{\varepsilon}))$. Hence, the argument used in the proof of the upper bound $\HD_d(p,q) \leq \tilde{O}(p^{\frac{q-1}{q-d}\cdot \rho_d})$ in~\cite{KST18} yields an upper bound of $\tilde{O}(p^{\frac{q-1}{q-d}})$ for such families. 

We note that both previous strategies for obtaining lower bounds on $\HD_d(p,d+1)$ cannot yield stronger lower bounds than the $p^{d-o(1)}$ bound of Theorem~\ref{thm:main}. Indeed, as was mentioned above, the approach via weak $\varepsilon$-nets cannot yield a bound stronger than $p^{\rho_d}$, for $\rho_d \approx d-0.5$. The approach via hyperplanes produces families that have a bounded VC dimension, and as explained above, the bound of Theorem~\ref{thm:main} is sharp within this class. 

\medskip \noindent \textbf{Open problem.} We do not know whether the lower bound of Theorem~\ref{thm:main} is tight. A natural approach to prove tightness in $\mathbb{R}^2$ is to prove that every finite point set in the plane admits a weak $\varepsilon$-net of size $\tilde O(1/\varepsilon)$. For $d \geq 3$, such an approach cannot yield tightness due to the aforementioned superlinear lower bound of Bhore, Browne, and Zheng~\cite{BBZ26} for the $\varepsilon$-net theorem. Hence, if the bound of Theorem~\ref{thm:main} is tight for $d \geq 3$, proving this will probably require devising a conceptually different proof of the $(p,q)$ theorem, a goal that is certainly interesting in itself. 

\subsection{Organization of the paper}

In Section~\ref{sec:overview} we present a detailed overview of the proof. In Section~\ref{sec:HD2(p,3)} we present the proof of the lower bound on $\HD_2(p,3)$. In Section~\ref{sec:HDd(p,d+1)} we present the generalization to $\HD_d(p,d+1)$ for any $d \geq 2$. Finally, in Section~\ref{sec:general} we present the lower bound on $\HD_d(p,q)$ for $q \geq d+2$.

\section{Proof overview}
\label{sec:overview}

We prove the bound by constructing a specific family of parallelotopes and using the \emph{container lemma}~\cite{BMS15,ST15} to show that it satisfies the $(p,q)$ property and its piercing number is at least the bound stated in the theorem. For the sake of clarity, we first describe in detail the construction and proof of the most basic case $d=2, q=3$, and then we briefly describe the way to extend it to $d>2$ and to $q>d+1$. 

\begin{additions}
\subsection{The lower bound on $\HD_2(p,3)$}
\end{additions}

\noindent \textbf{A random family of strips.} Our first idea is to replace the families of lines used in the earlier constructions in the plane (\cite{CGS26,KellerS21,RN26}) with a random family of thin strips. For $i=1,\ldots,n$, we sample $\{P_i=(a_i,b_i)\}_{i=1,2,\ldots,n}$ independently and uniformly from $[0,1]^2$ and set the following family of strips:
\[
  S_i=\{(x,y)\in\R^2: |y-a_ix-b_i|\le 1/n\},
  \qquad \mathcal S_n=\{S_1,\ldots,S_n\}.
\]
We will show that with high probability, the intersections of these strips with a sufficiently large square forms a family of compact parallelograms that proves the assertion of Theorem~\ref{thm:main}. 

We call the unit square $[0,1]^2$, where we sample the coefficient pairs $(a_i,b_i)$,  the \emph{coefficient space}. Note that while the coefficients are taken only in the unit square, the sets $S_i$ are thin bands in the whole $(x,y)$-plane. To prevent confusion, we call these bands \emph{strips}, and use the term \emph{tubes} for bands about a line in the coefficient space.

\medskip \noindent \textbf{The general structure of the proof.} We prove that with probability tending to one as $n\to\infty$, the family $\mathcal S_n$ satisfies the two following properties:
\begin{enumerate}
\item[(A)] Every point of the plane belongs to fewer than
  $r_n=\lceil 8\log n/\log\log n\rceil$ strips.
\item[(B)] Every subfamily of $p_n=\lceil\sqrt n(\log n)^3\rceil$ strips contains three strips with a common point.
\end{enumerate}
By (A), piercing all the strips requires at least $n/r_n$ points. For a given $p$, choosing $n$ of order $p^2/(\log p)^6$ with a sufficiently small constant makes $p_n\le p$, and thus, (B) gives the $(p,3)$ property. The number of points required to pierce the family, which is at least $n/r_n$, has the order of the lower bound stated in the theorem.

\medskip \noindent \textbf{Proving (A) using a geometric observation and a union bound.} Proving (A) is the easier part. We observe that if some point belongs to $t$ strips, then there exists a point in the intersection of two boundaries of strips that also belongs to $t$ strips (see Figure~\ref{fig:depth} below). The number of such intersection points is at most $4\binom{n}{2}$, and for each such point, we can bound from above the probability that it belongs to at least $t-2$ additional strips using the independence between strips and the fact that each other strip contains the point with probability about $2/n$. The assertion $t \leq r_n$ follows from a union bound over the $4\binom{n}{2}$ intersection points.  

\medskip \noindent \textbf{Transition from intersections in the plane to co-linearity in the coefficient space.} The basic observation behind the proof of (B) is that intersections between strips in the plane can be translated to almost co-linearity of the corresponding pairs of coefficients $(a_i,b_i)$ in the coefficient space. We observe that if three coefficient points $\{(a_i,b_i)\}_{i=1,2,3}$ lie on the line $b=ua+v$, then the central lines $y=a_ix+b_i$ of the corresponding strips all pass through $(-u,v)$. The three strips still contain $(-u,v)$ whenever their coefficient points lie within vertical distance $1/n$ of the line $b=ua+v$. Hence, a \emph{tube} of vertical halfwidth $1/n$ about the line $b=ua+v$ in the coefficient space, which we denote by $T_{u,v}(1/n)$, consists of the coefficients of strips that contain the point $(-u,v)$ in the plane, as illustrated in Figure~\ref{fig:duality}.

\medskip \noindent \textbf{Reduction to independent sets in a $3$-uniform hypergraph.}
For the proof of (B), we divide the coefficient square $[0,1]^n$ into $V=n^6$ squares of side length $n^{-3}$, which we call \emph{cells}. Suppose a tube $T_{u,v}(1/n)$ contains three entire cells. Then any choice of one coefficient point from each of the three cells gives three strips that all contain the point $(-u,v)$.

Before the coefficient points $P_i$ are sampled, we define $H$ to be the $3$-uniform hypergraph whose vertices are the cells. Three distinct cells form an edge of $H$ if they are contained in a single tube of the form $T_{u,v}(1/n)$. A subfamily of strips with no intersecting triple corresponds to sampled coefficient points whose cells form an \emph{independent set} of $H$, namely, a set of vertices of $H$ containing no edge. (We use here the fact that with high probability, the sampled points come from distinct cells.) Therefore, in order to prove (B), it  suffices to show that the size of any independent set in $H$ that is contained in our sample is less than $p_n$.

\medskip \noindent \textbf{Using the hypergraph container lemma.} Since its introduction in 2012, the hypergraph container lemma~\cite{BMS15,ST15} has become a central method for bounding the size of independent sets in hypergraphs. In our context, it can be applied to a $3$-uniform hypergraph if it is sufficiently dense and its edges are `evenly spread' among the vertices, in the sense that no vertex and no pair of vertices are included in `too many' edges. The lemma allows constructing a `not very large' family $\mathcal{C}$ of sets of vertices called \emph{containers}, such that each container is `not too large', and every independent set in the hypergraph is contained in some container $C \in \mathcal{C}$.

We use a sequence of iterative applications of the container lemma to $H$ and its subhypergraphs (see below) to obtain a family $\mathcal{C}$ of at most $\exp(\sqrt n\,(\log n)^{5/2})$ containers in $H$, each of size at most $\sqrt{\log n/n}\cdot V$, such that each independent set is contained in one of them. We do this \emph{before} sampling the points $\{(a_i,b_i)\}_{i=1,2,\ldots,n}$. 

Each container $C$ occupies area $|C|/V$ in the coefficient square, so the expected number of sample points whose cells are included in it is about $n|C|/V$. By the upper bound on the size of the container, this expectation is much smaller than $p_n$. Hence, the probability that the container contains at least $p_n$ points is tiny, and can be bounded using tail estimates for binomial random variables. As the number of containers is not too large, we can use a union bound to show that with high probability, no container contains at least $p_n$ sample points. Hence, the cells corresponding to the sampled points in $[0,1]^2$ contain no independent set of size $p_n$, which implies that the corresponding strips satisfy the $(p_n,3)$-property, proving (B).

\medskip \noindent \textbf{Iterative application of the container lemma to construct $\mathcal{C}$.}
The construction process of $\mathcal C$ is iterative. It starts with a single container of all the cells and repeatedly replaces each large container $J$ with a family of smaller containers included in it, which we call its \emph{children}. Each such replacement is an application of the container lemma to a specifically chosen subhypergraph $H_J \subset H$ whose vertex set is the respective $J$. The form of the container lemma we use (Lemma~\ref{lem:container} below) ensures that every independent set of its input hypergraph is included in some child and that each child retains at most $1/4$ of that hypergraph's edges. One of the properties of $H_J$ is that each cell has degree at most a constant times the average degree of $H_J$ (see below). Hence, the edge reduction guaranteed by the container lemma forces each child to omit at least a fixed fraction of the cells of $J$. Since $H_J$ is a subhypergraph of $H[J]$, every independent set of $H$ contained in $J$ is independent in $H_J$ and is therefore contained in a child. We repeat the step until the containers reach the desired size.

\medskip \noindent \textbf{Application of the container lemma to each subhypergraph $H_J$.} The most complex and technical part of the proof is choosing the $3$-uniform subhypergraphs $H_J$ and proving that each $H_J$  satisfies the conditions required for the application of the container lemma. Specifically, three conditions must be satisfied:
\begin{enumerate}
\item[(i)] $|E(H_J)| = \Omega(|J|^3/n \log n)$;
\item[(ii)] The degree of each cell in $H_J$ is at most constant times the average vertex degree of $H_J$;
\item[(iii)] Each pair of cells is contained in at most $O(V/n)$ edges of $H_J$.
\end{enumerate}
We construct $H_J$ edge-by-edge according to a procedure described below. To make sure (ii) is satisfied, we fix $D$, maintain a subset $J'$ of all cells whose total degree hasn't reached $D$, and update the subset after each step. By choosing $D$ properly and adding only edges all whose cells are included in $J'$ at the respective step, we guarantee that (ii) is satisfied. To make sure (iii) is satisfied, we check at each step that all three pairs of cells within the triple we want to add are contained in at most $O(V/n)$ edges chosen so far. However, how can we guarantee that the process can be continued for $\Omega(|J|^3/n \log n)$ steps and that at each step, there are available triples satisfying the conditions we check? 

Here, the procedure of adding edges comes into play. The key is handling condition (iii). We observe that for any fixed $w>0$, pairs of points in the \emph{coefficient space} whose $x$ coordinates are close to each other, belong to much more tubes of the form $T_{u,v}(w)$ than pairs of points whose $x$ coordinates are far from each other (see Figure~\ref{fig:gaps}). Specifically, if the union of the cells in $J$ is $A$, then the tubes containing a pair of points with horizontal gap $g$ occupy a parameter region of area $4w^2/g$. To compensate for this, for each tube, we write $2wL$ for the area of its intersection with $A$,
and include in the hypergraph $H_J$ \emph{only triples of `far' points}, specifically, triples whose consecutive horizontal gaps are at least $L/8$. Thus, a tube that contributes to the number of edges in $H_J$ that contain a pair with gap $g$, must have $L\le 8g$. The area available for the third point is then at most $16wg$, which compensates for the factor $1/g$ in the area of the tubes that contain the two points.

Based on this observation and on averaging over all points in a cell, for each $J' \subset J$ we define a somewhat complex \emph{weight function} $\eta(J')$ on triples of cells in $J'$ which makes sure that the overall weight of the triples of cells that contain a pair of cells is spread close-to-uniformly over the pairs of cells. Then, we construct the hypergraph $H_J$ edge by edge, where at each step, we use the weight function $\eta(J')$ to sample an edge of $H[J']$ and add it to $H_J$. The weight function makes it possible to prove that at each step, there are available triples satisfying the conditions we check and that the process can be continued until sufficiently many edges are collected. The latter proof is rather technical and uses calculations of tube areas.

\medskip

\begin{additions}
\subsection{The lower bound on $\HD_d(p,d+1)$}

\textbf{The general structure of the proof.} The general structure is the same as in the planar case. We sample coefficient points $\{P_i=(a_i,b_i)\}_{i=1,\ldots,n}$, now in $[0,1]^d$,
with each $a_i \in [0,1]^{d-1}$ a vector of $d-1$ slopes. The strips become thin slabs, and $T_{u,v}(1/n)$ is the tube about the hyperplane $b=u\cdot a+v$ in the coefficient space. As in the plane, if the tube contains $d+1$ entire cells, every choice of one point in each cell gives $d+1$ slabs that contain the point $(-u,v) \in [0,1]^d$. Hence, we can define a hypergraph $H$ whose vertices are the cells and whose hyperedges are $(d+1)$-tuples of cells contained in a tube of the form $T_{u,v}(1/n)$, and use the container lemma to bound the maximum size of an independent set in $H$ that is included in our sample, and thus to prove the analogue of (B).

\medskip \noindent \textbf{A more complex weight function.} The main difficulty in the generalization to all $d \geq 2$ is that the `uniform-spreadness' condition (iii) above extends from pairs of cells to every fixed set of two to $d$ cells of $H_J$, and is generally harder to handle in $\mathbb{R}^d$. In the plane, a large horizontal gap between two coefficient points implies that the number of tubes $T_{w}(u,v)$ that contain the two points is small. In dimension $d=3$, triples of points whose slope vectors are well-separated can be included in many tubes. An example is the slope vectors $(0,0), (1,0), (1/2, \varepsilon)$, and the containment in many tubes occurs since the point $(1/2,\varepsilon)$ is close to the line connecting $(0,0),(1,0)$. More precisely, for $d=3$, the measure of the set of tubes containing the three points is inversely proportional to the area of the triangle formed by the three slope vectors. For $d$ coefficient points with affinely independent slope vectors in $\mathbb{R}^d$, the corresponding measure is
inversely proportional to the $(d-1)$-dimensional volume of their slope simplex. If we fix only $j<d$ coefficient points whose slope vectors span a $(j-1)$-dimensional affine subspace, the center of the tube can still tilt in the $d-j$ directions perpendicular to this subspace. When we add $d-j$ points to complete the simplex, its volume is, up to a dimensional constant, the volume of the fixed face times the volume spanned by the perpendicular components of
the added slope vectors. 

We use this factorization to bound the total weight of completions of the fixed $j$ points. We set the weight of a $(d+1)$-tuple to combine the measure of the tubes that contain it with the largest $(d-1)$-dimensional simplex volume among its slope vectors. Expanding the determinant that gives this volume in the directions
perpendicular to the fixed face produces terms that cancel the corresponding inverse factors in the tube measure. The remaining inverse factor is the volume of the fixed face, so we impose a lower bound on the volume of every face. If the slope vectors concentrate near a flat, we first rescale the narrow directions
to ensure that enough tuples pass these volume cutoffs. We then average the weights over cells and use the bounds on the total weight and on the weight of completions to select many edges of $H_J$ with the required degree bounds.
\end{additions}

\subsection{The lower bound on $\HD_d(p,q)$}

The construction extends to every fixed $q\ge d+2$ using the same random slabs. We now use a $q$-uniform hypergraph in which an edge consists of $q$ cells contained in a common tube of the form $T_w(u,v)$. The weight of a tuple still uses the largest simplex volume among its $d$-point faces. The $q-d$ additional points must all lie in the same tube, so the factor $w$ in the total weight becomes $w^{q-d}$. A
small-determinant estimate shows that sufficiently many choices of the additional points satisfy the volume cutoffs on all faces simultaneously. To bound the weight of completions of $j$ fixed points, we use the same determinant cancellation for $2\le j\le d$.
When $d<j<q$, a $d$-point face among the fixed points already bounds the measure of the tubes containing all points, and integrating the remaining $q-j$ points within each tube contributes a factor $w^{q-j}$.
The resulting container bounds show that, with high probability, the largest subfamily containing no intersecting $q$-tuple has size
$O_{d,q}(n^{(q-d)/(q-1)}(\log n)^3)$. Together with the depth bound (i.e., the analogue of (A) above),
this gives $\HD_d(p,q)\ge p^{(q-1)/(q-d)-o(1)}$.

\section{The lower bound on $\HD_2(p,3)$}
\label{sec:HD2(p,3)}

In this section we prove the lower bound for the $(p,3)$ theorem in the plane. An overview of the proof was given in Section~\ref{sec:overview} above.  
Section~\ref{sec:random} proves assertion (A) from Section~\ref{sec:overview} and shows how (A) and (B) yield the assertion of the theorem. Section~\ref{sec:cells} defines the cell hypergraph $H$ and deduces (B) from the container bounds. Section~\ref{sec:geometric} constructs the auxiliary hypergraphs $H_J$ and proves the degree bounds. Section~\ref{sec:containers} uses these bounds to construct the containers and complete the proof of (B). 

\Needspace{9\baselineskip}
\subsection{Random strips and the piercing bound}\label{sec:random}

Almost surely, any intersecting subfamily of at least two strips has
a common point at the intersection of two of its boundary lines.
These boundary intersections give a finite list of candidates for the
union bound in (A) from Section~\ref{sec:overview}. They also let us choose a single square that
preserves every intersecting triple when we truncate the strips.

For $1\le i\le n$, choose $P_i=(a_i,b_i)$ independently and uniformly from $[0,1]^2$, and set
\[
  S_i=\{(x,y)\in\R^2: |y-a_ix-b_i|\le 1/n\},
  \qquad \mathcal S_n=\{S_1,\ldots,S_n\}.
\]
The halfwidth $1/n$ is measured vertically. For a finite family $\mathcal F$, write $\alpha_3(\mathcal F)$ for the largest size of a subfamily containing no three distinct members with a common point, and write
\[
  \depth(\mathcal F)=\max_{z\in\R^2}|\{F\in\mathcal F:z\in F\}|.
\]
If $|\mathcal F|\ge p$, the $(p,3)$ property is equivalent to $\alpha_3(\mathcal F)<p$. Also, each piercing point meets at most $\depth(\mathcal F)$ members. We will use the thresholds
\[
  p_n=\left\lceil\sqrt n\,(\log n)^3\right\rceil,
  \qquad
  r_n=\left\lceil\frac{8\log n}{\log\log n}\right\rceil.
\]
All logarithms are natural. In this notation, (A) is $\depth(\mathcal S_n)<r_n$ and (B) is
$\alpha_3(\mathcal S_n)<p_n$. We prove (A) below and use (B) in the
truncation argument. The proof of (B) is given in
Sections~\ref{sec:cells}--\ref{sec:containers}.

\paragraph{Depth.}
A fixed point $(x,y)$ lies in a random strip with probability at most $2/n$. Indeed, after conditioning on its slope $a$, the allowable intercepts lie in an interval of length $2/n$. To bound the depth throughout the plane, we use the intersections of strip boundaries as candidate points.

\begin{figure}[htbp]
\centering
\begin{tikzpicture}[x=1.05cm,y=1.05cm,font=\small]
  \begin{scope}
    \clip (-3.7,-2.43) rectangle (3.7,2.43);
    \fill[sea!6] (-3.7,-.50) -- (3.7,.24) -- (3.7,.66) -- (-3.7,-.08) -- cycle;
    \fill[indigo!7] (-3.7,-1.484) -- (3.7,1.624) -- (3.7,2.044) -- (-3.7,-1.064) -- cycle;
    \fill[ochre!8] (-3.7,-2.58) -- (3.7,2.60) -- (3.7,3.02) -- (-3.7,-2.16) -- cycle;
    \draw[strip a,opacity=.43] (-3.7,-.50) -- (3.7,.24)
                                 (-3.7,-.08) -- (3.7,.66);
    \draw[strip b,opacity=.45] (-3.7,-1.484) -- (3.7,1.624)
                                 (-3.7,-1.064) -- (3.7,2.044);
    \draw[strip c,opacity=.46] (-3.7,-2.58) -- (3.7,2.60)
                                 (-3.7,-2.16) -- (3.7,3.02);
    \begin{scope}
      \clip (-2.3,-2.3) rectangle (2.3,2.3);
      \fill[sea!13] (-3.7,-.50) -- (3.7,.24) -- (3.7,.66) -- (-3.7,-.08) -- cycle;
      \fill[indigo!14,fill opacity=.8] (-3.7,-1.484) -- (3.7,1.624) -- (3.7,2.044) -- (-3.7,-1.064) -- cycle;
      \fill[ochre!18,fill opacity=.76] (-3.7,-2.58) -- (3.7,2.60) -- (3.7,3.02) -- (-3.7,-2.16) -- cycle;
      \fill[ink!28] (-.93333,-.22333) -- (-.23333,.26667) -- (.46667,.33667)
          -- (.21429,.16) -- (-.625,-.1925) -- cycle;
      \draw[strip a] (-3.7,-.50) -- (3.7,.24)
                     (-3.7,-.08) -- (3.7,.66);
      \draw[strip b] (-3.7,-1.484) -- (3.7,1.624)
                     (-3.7,-1.064) -- (3.7,2.044);
      \draw[strip c] (-3.7,-2.58) -- (3.7,2.60)
                     (-3.7,-2.16) -- (3.7,3.02);
    \end{scope}
  \end{scope}
  \draw[ink,thick] (-2.3,-2.3) rectangle (2.3,2.3);
  \node[ink,below right] at (-2.3,2.3) {$B$};
  \fill[ink] (.46667,.33667) circle[radius=1.7pt];
  \draw[ink,semithick,-{Latex[length=1.7mm]}] (.88,1.23) -- (.492,.405);
  \node[ink,above,inner sep=1pt] at (.88,1.23) {$z$};
\end{tikzpicture}
\caption{An extreme point $z$ of the common intersection lies on two strip boundaries. After the sample is fixed, the square $B$ is chosen to contain all intersections of nonparallel boundaries. In the drawing each strip is cut only by the vertical sides of $B$.}
\label{fig:depth}
\end{figure}
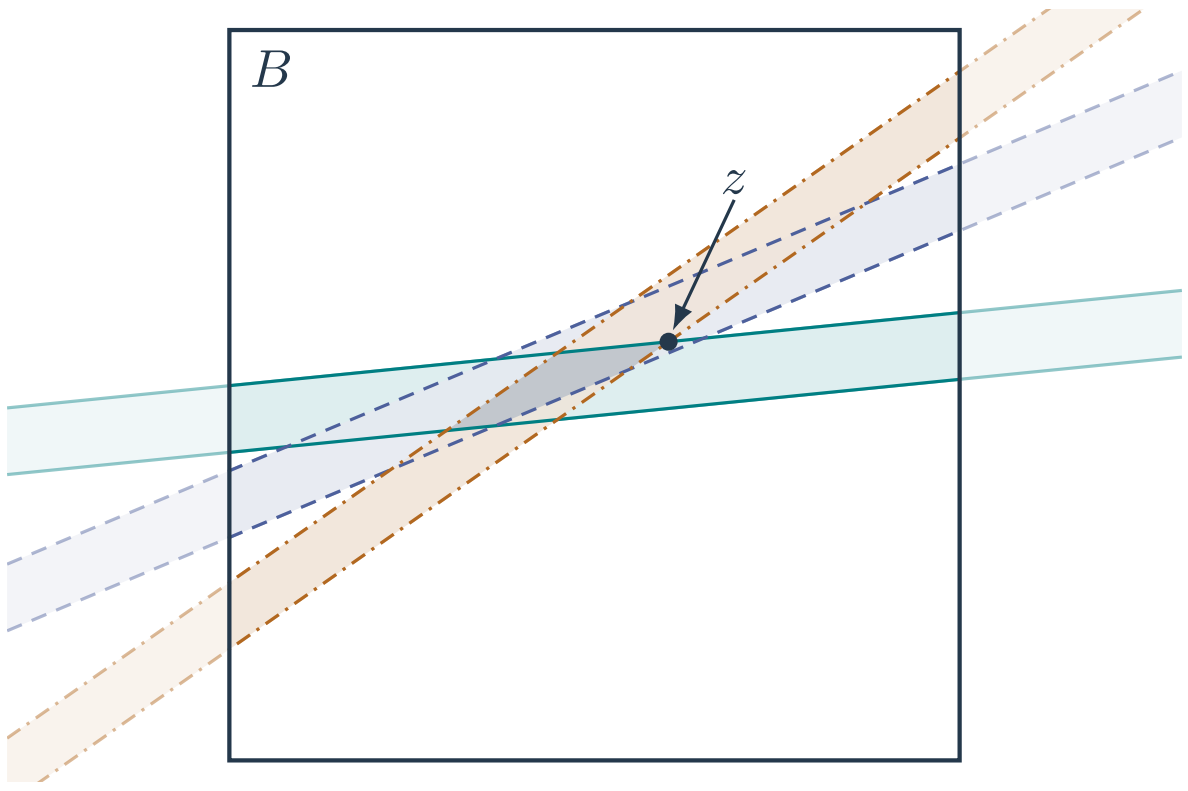

The slopes are distinct almost surely. If a subfamily of at least two strips has a common point, its intersection is compact because any two of the strips intersect in a compact parallelogram. At an extreme point of the intersection, two active boundary lines are nonparallel. Otherwise one could move slightly in both directions along the common direction of the active lines and stay in the intersection. Thus any such subfamily has a common point at the intersection of boundaries of two of its members (Figure~\ref{fig:depth}).

There are at most $4\binom n2$ intersections of boundaries from two different strips. Conditional on the two strips defining one candidate, the remaining strips are independent, each containing that point with probability at most $2/n$. For $3\le r\le n$ the union bound gives
\[
  \Prob\bigl(\depth(\mathcal S_n)\ge r\bigr)
  \le 4\binom n2\binom{n-2}{r-2}(2/n)^{r-2}
  \le 2n^2\left(\frac{2e}{r-2}\right)^{r-2}.
\]
At $r=r_n$ the last expression tends to zero, since $(r_n-2)\log(r_n-2)=(8+o(1))\log n$.

\paragraph{Compactness and the choice of $n$.}
Assume for the moment that $\alpha_3(\mathcal S_n)<p_n$ with probability tending to one. Fix a sample satisfying both estimates, with distinct slopes and intercepts and $\max_i a_i<1$. Let $B=[-M,M]^2$ contain all intersections of nonparallel strip boundaries, with $M$ large enough that
\[
  \bigl(1-\max_i a_i\bigr)M>1+1/n .
\]
For $|x|\le M$, every $(x,y)\in S_i$ has $|y|\le a_i M+1+1/n<M$. Thus $F_i=S_i\cap B$ is a compact parallelogram with two vertical sides and nonempty interior. The $F_i$ are distinct since their traces on $x=0$ are distinct intervals. Every intersecting triple of strips has a common point at an intersection of two boundaries, so it still intersects after truncation. The truncated family $\mathcal F_n=\{F_1,\ldots,F_n\}$ therefore satisfies
\[
  \alpha_3(\mathcal F_n)<p_n,\qquad
  \tau(\mathcal F_n)\ge \frac{n}{r_n}
     =\Omega\!\left(\frac{n\log\log n}{\log n}\right).
\]
For a given sufficiently large $p$, take $n=\lfloor c_0p^2/(\log p)^6\rfloor$ with a sufficiently small absolute $c_0>0$. Then $n\ge p$, and $\log n\le 2\log p$ gives $\sqrt n(\log n)^3\le 8\sqrt{c_0}\,p$, so $p_n\le p$. As $n\asymp p^2/(\log p)^6$ and $\log n\asymp\log p$, the displayed bound for $\tau(\mathcal F_n)$ gives \added{the planar case of }Theorem~\ref{thm:main}. It remains to prove the estimate for $\alpha_3(\mathcal S_n)$.

\subsection{The cell hypergraph and the sampling step}\label{sec:cells}

To prove (B), we record which triples of cells guarantee a common point
of the corresponding strips for every choice of representatives. This
lets us define the hypergraph $H$ before sampling.

A point $z=(-u,v)$ lies in the strip associated with $(a,b)$ exactly when $|b-ua-v|\le 1/n$. Thus a common point of strips corresponds to a tube around a line in coefficient space. We write
\[
  T_{u,v}(w)=\{(a,b)\in\R^2: |b-ua-v|\le w\}
\]
for such a tube of vertical halfwidth $w$. The part of this tube in
coefficient space is $T_{u,v}(w)\cap[0,1]^2$.

\begin{figure}[htbp]
\centering
\begin{tikzpicture}[font=\small]
  \begin{scope}[x=4.75cm,y=4.12cm]
    \begin{scope}
      \clip (0,0) rectangle (1,1);
      \fill[sea!12] (0,-1.155) -- (1,1.845) -- (1,2.055) -- (0,-.945) -- cycle;
      \draw[step=.05,ink!10,very thin] (0,0) grid (1,1);
      \draw[sea!80!ink,semithick,dashed] (0,-1.05)--(1,1.95);
      \filldraw[fill=sea!33,draw=sea,thick] (.4,.2) rectangle (.45,.25);
      \filldraw[fill=indigo!27,draw=indigo,thick] (.5,.5) rectangle (.55,.55);
      \filldraw[fill=ochre!29,draw=ochre,thick] (.6,.8) rectangle (.65,.85);
    \end{scope}
    \draw[ink!60] (0,0) rectangle (1,1);
    \draw[axis] (0,0) -- (1.09,0) node[right] {$a$};
    \draw[axis] (0,0) -- (0,1.08) node[above] {$b$};
    \node[ink,below] at (1,0) {$1$};
    \node[ink,left] at (0,1) {$1$};
    \node[ink,anchor=east,inner sep=1pt] at (.37,.222) {$C_1$};
    \node[ink,anchor=east,inner sep=1pt] at (.48,.526) {$C_2$};
    \node[ink,anchor=east,inner sep=1pt] at (.58,.82) {$C_3$};
    \node[ink,anchor=west] at (.66,.3) {$T_{u,v}(1/n)$};
    \draw[ink!70,-{Latex[length=1.6mm]}] (.72,.37) -- (.567,.58);
    \node[ink] at (.5,-.15) {Coefficient space};
  \end{scope}
  \begin{scope}[shift={(10.03cm,2.50cm)},x=.89cm,y=1.29cm]
    \begin{scope}
      \clip (-4.52,-1.88) rectangle (1.57,1.26);
      \fill[sea!15] (-4.52,-1.801) -- (1.57,.78725) -- (1.57,.99725) -- (-4.52,-1.591) -- cycle;
      \fill[indigo!16,fill opacity=.7] (-4.52,-1.953) -- (1.57,1.24425) -- (1.57,1.45425) -- (-4.52,-1.743) -- cycle;
      \fill[ochre!22,fill opacity=.68] (-4.52,-2.105) -- (1.57,1.70125) -- (1.57,1.91125) -- (-4.52,-1.895) -- cycle;
      \draw[strip a] (-4.52,-1.801) -- (1.57,.78725)
                     (-4.52,-1.591) -- (1.57,.99725);
      \draw[strip b] (-4.52,-1.953) -- (1.57,1.24425)
                     (-4.52,-1.743) -- (1.57,1.45425);
      \draw[strip c] (-4.52,-2.105) -- (1.57,1.70125)
                     (-4.52,-1.895) -- (1.57,1.91125);
    \end{scope}
    \draw[axis] (-4.5,0) -- (1.64,0) node[right] {$x$};
    \draw[axis] (0,-1.8) -- (0,1.3) node[above] {$y$};
    \fill[ink] (-3,-1.05) circle[radius=1.6pt];
    \draw[ink,semithick,-{Latex[length=1.6mm]}] (-3.42,-.53) -- (-3.045,-1.005);
    \node[ink,above] at (-3.42,-.53) {$z=(-u,v)$};
    \node[ink] at (-1.47,-2.4186) {$(x,y)$-plane};
  \end{scope}
\end{tikzpicture}
\caption{The tube on the left contains three entire cells. The point
$z=(-u,v)$ belongs to all three strips for every choice of one
coefficient point from each cell. The strips on the right are drawn
using the cell centers. Widths and cell sizes are schematic.}
\label{fig:duality}
\end{figure}
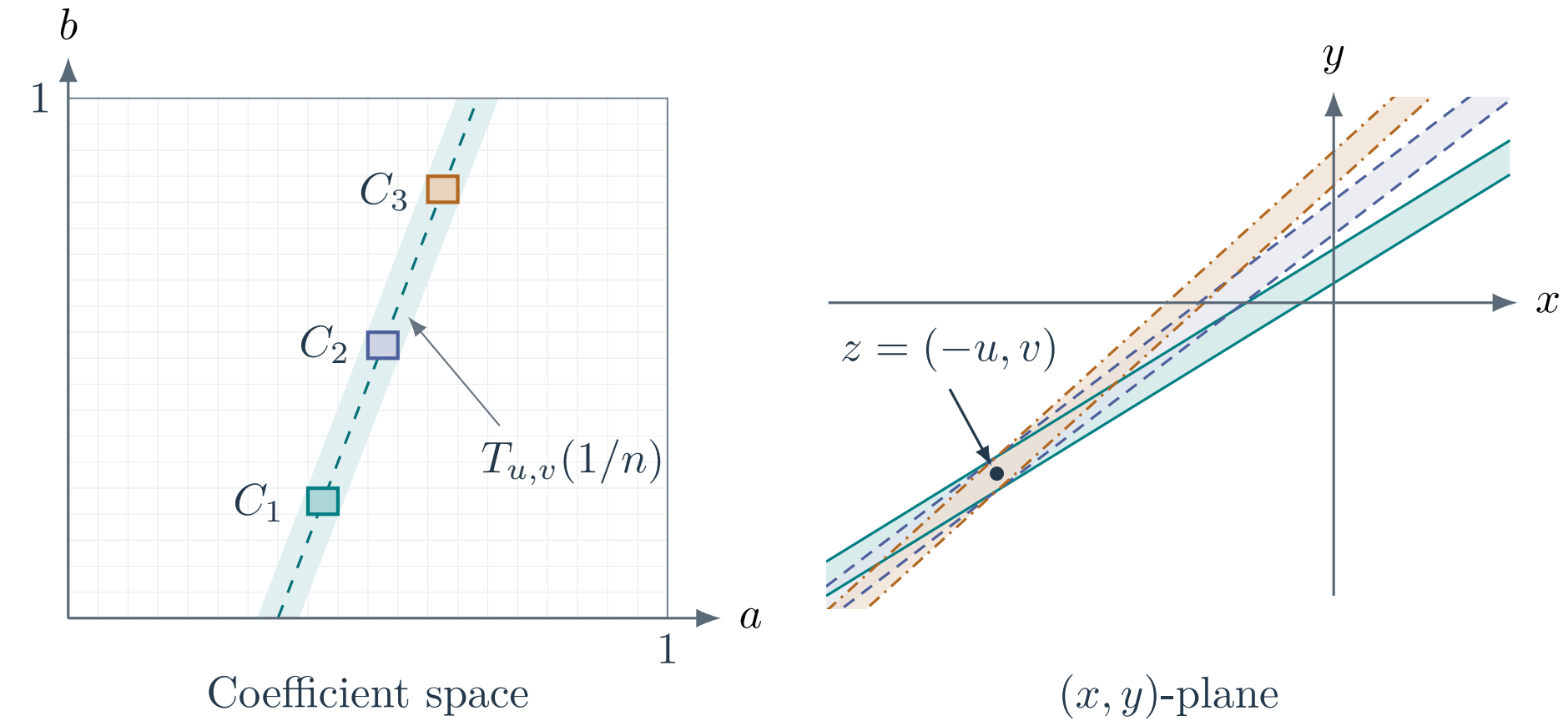

Partition $[0,1]^2$ into $V=n^6$ square cells of side $\eta=n^{-3}$, assigning each boundary point to just one cell. Let $H$ be the $3$-uniform hypergraph on the cells. Three distinct cells form an edge if some tube $T_{u,v}(1/n)$ contains all three cells. Equivalently, there is a single point $z=(-u,v)$ in the $(x,y)$-plane belonging to the strips for \emph{every} choice of one representative from each cell (Figure~\ref{fig:duality}). This definition makes $H$ independent of the sample.

\paragraph{From containers to (B).}
The area of a container controls how many sample points it is likely
to receive. The number of containers controls the cost of taking a
union bound. Proposition~\ref{prop:containers} gives both bounds.

\begin{proposition}\label{prop:containers}
For all sufficiently large $n$, there is a deterministic family $\mathcal C$ of containers covering the independent sets of $H$ such that
\[
  \max_{C\in\mathcal C}\frac{|C|}{V}
     =O\!\left(\sqrt{\frac{\log n}{n}}\right),
  \qquad
  \log|\mathcal C|=O\!\left(\sqrt n\,(\log n)^{5/2}\right).
\]
\end{proposition}

To see how this implies the desired estimate, fix the containers before drawing the $n$ coefficient points. The probability of any cell collision is at most $\binom n2/V=O(n^{-4})$. In the absence of collisions, any subfamily with no intersecting triple occupies distinct cells forming an independent set of $H$. Hence its size is at most the largest number of samples in any container.

For a fixed $C\in\mathcal C$, the number of sampled points in its cells has the binomial distribution with parameters $n$ and $|C|/V$, and mean at most $A\sqrt{n\log n}$ for an absolute $A$. This distribution is unconditional. By the binomial upper-tail bound and a union bound over $C$, the probability that some container holds at least $p_n$ samples is at most
\[
  |\mathcal C|
  \left(\frac{eA\sqrt{n\log n}}{p_n}\right)^{p_n}
  =o(1).
\]
The negative term in the logarithm has order $p_n\log\log n$, which dominates $\log|\mathcal C|$. Adding the collision probability proves $\Prob(\alpha_3(\mathcal S_n)\ge p_n)=o(1)$. The rest of the proof establishes Proposition~\ref{prop:containers}.

\subsection{Many triples with controlled degrees}\label{sec:geometric}

For a $3$-uniform hypergraph $G$, let $e(G)$ be its number of edges and let $\Delta_j(G)$ be the largest number of edges containing any fixed set of $j$ vertices. At each step of the container construction, $J$ is the current collection of cells. We select many edges of $H[J]$, with a bound on the number of selected edges through any pair of cells so that we can apply the container theorem. We also bound the number through any single cell, so that discarding a fixed fraction of the selected edges requires discarding a fixed fraction of the cells.

The bounds on $\Delta_1(H_J)$ and $\Delta_2(H_J)$ below correspond
to (i) and (ii) of the introduction, respectively.

\Needspace{7\baselineskip}
\begin{proposition}\label{prop:balanced}
There are absolute constants $c,C>0$ such that the following holds for all sufficiently large $n$. For every collection $J$ of $m\ge 2V/\sqrt n$ cells, there is a subhypergraph $H_J$ of $H[J]$, on the vertex set $J$, with
\[
  e(H_J)\ge c\,\frac{m^3}{n\log n},
  \qquad
  \Delta_1(H_J)\le C\,\frac{e(H_J)}{m},
  \qquad
  \Delta_2(H_J)\le C\,\frac{V}{n}.
\]
\end{proposition}

\paragraph{Idea of the construction.}
In Section~\ref{subsec:weighted}, we work with the union $A$ of the cells of $J$ and assign a weight $K$ to triples of points. A pair with horizontal gap $g$ lies in tubes $T_{u,v}(w)$ whose parameters $(u,v)$ fill a region of area $4w^2/g$. Close pairs thus have more possible witnesses. We count a triple in a tube only if its consecutive gaps are at least $L/8$, where $2wL$ is the area of $A$ inside the tube. A tube that contributes to a pair of gap $g$ must then have $L\le 8g$, so the area available for the third point is at most $16wg$. The factors $1/g$ and $g$ cancel (Figure~\ref{fig:gaps}). Lemma~\ref{lem:weighted} shows that this rule keeps the weight through each pair small while retaining a large total weight.

In Section~\ref{subsec:rounding}, we take $w=1/(100n)$. The slope bound in Lemma~\ref{lem:weighted} ensures that widening a witnessing tube from halfwidth $w$ to $1/n$ covers the three entire cells. We average the point weights over the cells to get a weight $\omega_E$ on each cell triple $E$, and keep each edge with probability $\omega_E$. The bound on the weight through a pair of points gives a bound on the expected degree of a pair of cells. Lemma~\ref{lem:round} gives many candidate edges with bounded pair degrees on every sufficiently large subcollection of $J$.

In Section~\ref{subsec:vertex-degrees}, we build $H_J$ one edge at a time. At each step we obtain candidates from Lemma~\ref{lem:round} on the cells whose degrees are still below their threshold. We need an unused candidate that avoids pairs already at their threshold. Each edge we add increases only three pair degrees, so few pairs reach their threshold. Each of those pairs occurs in at most $2V/n$ candidates, leaving an eligible edge at every step until $H_J$ reaches the desired size.

\subsection{A weighted estimate for points in a square}\label{subsec:weighted}

Let $A\subseteq[0,1]^2$ be a Borel set of area $\alpha>0$. In the application, $A$ is the union of the cells of $J$. Write $P=(a(P),b(P))$ for the coordinates of a point. Integrals over points and tube parameters are with respect to Lebesgue measure.

We assign weight to a triple by integrating over the tubes that witness
it. To bound the total weight from below, we fix a tube and count the
triples it contributes. To bound the weight of the triples containing
a given pair, we fix that pair and estimate how much each of its
witnessing tubes can contribute through a third point.

\begin{lemma}\label{lem:weighted}
For $0<w\le 1/10$ there is a symmetric measurable weight $K$ on triples of points in $[0,1]^2$, taking values in $[0,1/2]$, with the following properties.
\begin{enumerate}
\item If $K(P_1,P_2,P_3)>0$, the three points lie in a common tube $T_{u,v}(w)$ with $|u|\le C/\alpha$. When the points are ordered by first coordinate, both consecutive gaps are at least $c\alpha/\log(2/\alpha)$.
\item The total weight on triples in $A$ ordered by first coordinate satisfies
\[
  \int_{\substack{P_1,P_2,P_3\in A\\ a(P_1)<a(P_2)<a(P_3)}}
       K(P_1,P_2,P_3)\diff P_1\diff P_2\diff P_3
  \ge c\,\frac{w\alpha^3}{\log(2/\alpha)}.
\]
\item Every pair $P,Q\in[0,1]^2$ satisfies $\displaystyle\int_A K(P,Q,Z)\diff Z\le 8w$.
\end{enumerate}
Here $c,C>0$ are absolute constants.
\end{lemma}

\begin{proof}
For a fixed tube $T_{u,v}(w)$, let
\[
  f(a)=\frac{1}{2w}\int_{-w}^{w}\one_A(a,ua+v+t)\diff t,
  \qquad L=\int_0^1 f(a)\diff a .
\]
We extend $\one_A$ by zero outside the square. Thus $0\le f\le 1$ and $2wL=\area(A\cap T_{u,v}(w))$. The section quantities $f$ and $L$ depend on $(u,v)$.
At horizontal coordinate $a$, the tube has a vertical segment of length
$2w$, and $f(a)$ is the fraction of that segment belonging to $A$.
In particular, an interval of horizontal coordinates of length $s$
contributes at most $s$ to $L$.

\paragraph{The separation rule.}
Fix two points $P,Q$ with horizontal gap $g=|a(P)-a(Q)|>0$. Let $\mathcal T(P,Q)$ be the set of tube parameters $(u,v)$ for which the tube contains both points. The map
\[
  (u,v)\longmapsto
  \bigl(b(P)-ua(P)-v,\ b(Q)-ua(Q)-v\bigr)
\]
has absolute determinant $g$. Both residuals must lie in $[-w,w]$, so $\area(\mathcal T(P,Q))=4w^2/g$. Thus a small horizontal gap gives a large region of possible witnesses.

For each tube, we count a triple only if its consecutive horizontal
gaps are at least $L/8$. The gap of any pair in the triple is either
one of the consecutive gaps or their sum. Thus if a retained triple
contains $P,Q$, then $g\ge L/8$. In each tube that can contribute to
such triples, the third point must lie in the section
$A\cap T_{u,v}(w)$, of area $2wL\le 16wg$. This factor of $g$ in the
completion area will cancel the factor of $1/g$ in the area of
$\mathcal T(P,Q)$. Figure~\ref{fig:gaps} shows the two regions.

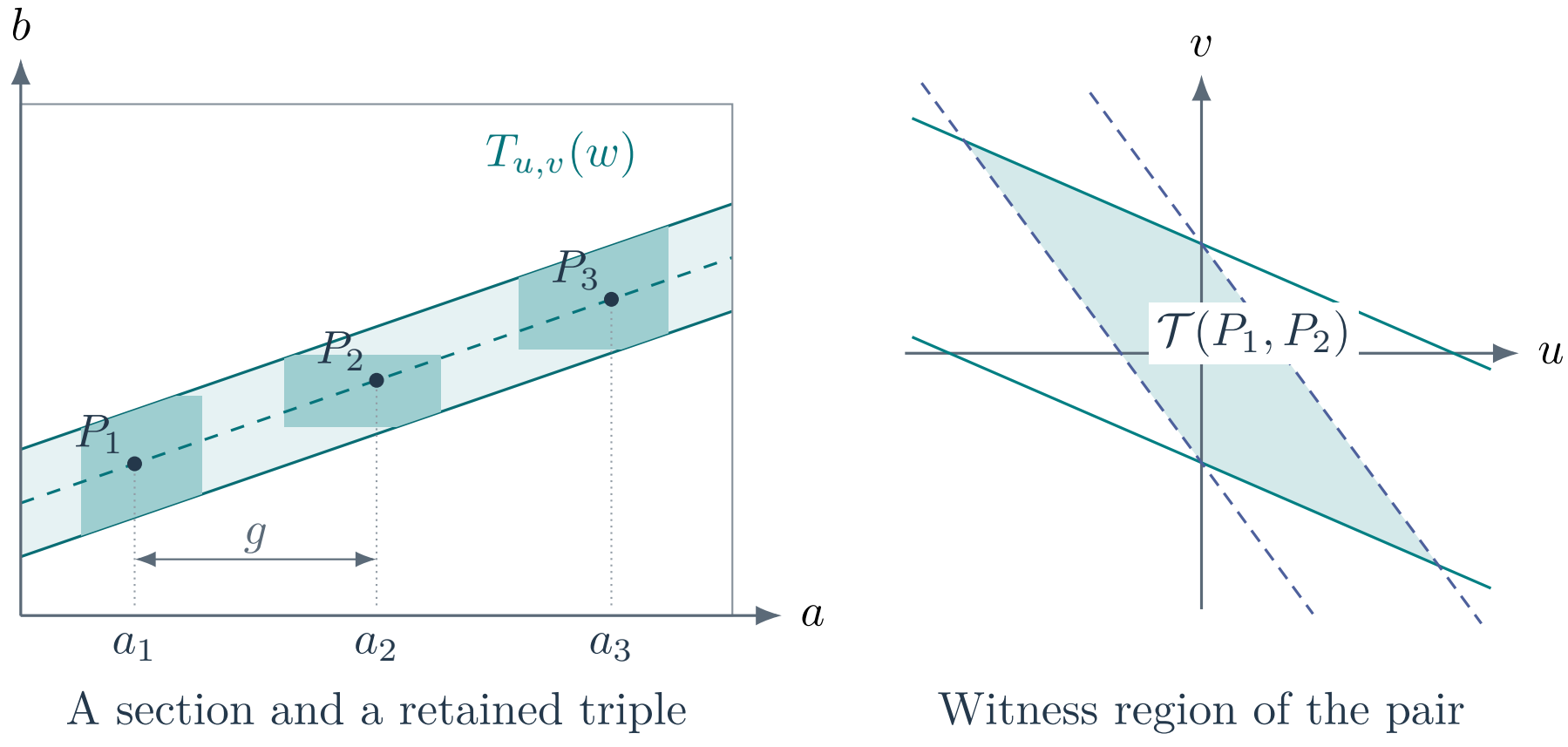
\begin{figure}[htbp]
\centering
\begin{tikzpicture}[font=\small]
  \begin{scope}[x=5.4cm,y=3.88cm]
    \fill[sea!10] (0,.115) -- (1,.595) -- (1,.805) -- (0,.325) -- cycle;
    \draw[sea!74!ink,semithick] (0,.115) -- (1,.595) (0,.325) -- (1,.805);
    \begin{scope}
      \clip (0,.115) -- (1,.595) -- (1,.805) -- (0,.325) -- cycle;
      \fill[sea!38] (.085,.15) rectangle (.255,.43);
      \fill[sea!38] (.37,.37) rectangle (.59,.51);
      \fill[sea!38] (.70,.52) rectangle (.91,.77);
    \end{scope}
    \draw[sea!85!ink,semithick,dashed] (0,.22) -- (1,.70);
    \draw[ink!55] (0,0) rectangle (1,1);
    \draw[axis] (0,0) -- (1.07,0) node[right] {$a$};
    \draw[axis] (0,0) -- (0,1.09) node[above] {$b$};
    \foreach \aa/\bb/\ii in {.16/.2968/1,.50/.46/2,.83/.6184/3} {
      \draw[ink!50,densely dotted] (\aa,.02) -- (\aa,\bb);
      \fill[ink] (\aa,\bb) circle[radius=1.6pt];
      \node[ink,above left,inner sep=2pt] at (\aa,\bb) {$P_{\ii}$};
      \node[ink,below] at (\aa,0) {$a_{\ii}$};
    }
    \draw[ink!75,{Latex[length=1.6mm]}-{Latex[length=1.6mm]}]
      (.16,.11) -- node[above,inner sep=1.2pt] {$g$} (.50,.11);
    \node[sea!85!ink] at (.76,.9) {$T_{u,v}(w)$};
    \node[ink] at (.5,-.19) {A section and a retained triple};
  \end{scope}
  \begin{scope}[shift={(8.96cm,1.99cm)},x=1.8cm,y=.83cm]
    \fill[sea!17] (0,-1) -- (1,-1.9412) -- (0,1) -- (-1,1.9412) -- cycle;
    \draw[axis] (-1.25,0) -- (1.34,0) node[right] {$u$};
    \draw[axis] (0,-2.34) -- (0,2.55) node[above] {$v$};
    \draw[strip a] (-1.22,.1482) -- (1.22,-2.1482)
                   (-1.22,2.1482) -- (1.22,-.1482);
    \draw[strip b] (-.47,2.3824) -- (1.18,-2.4706)
                   (-1.18,2.4706) -- (.47,-2.3824);
    \node[ink,fill=white,inner sep=1.6pt] at (.22,.18) {$\mathcal T(P_1,P_2)$};
    \node[ink] at (0,-3.2846) {Witness region of the pair};
  \end{scope}
\end{tikzpicture}
\caption{The darker portions of the tube on the left have total area $2wL$. Both consecutive gaps of a retained triple are at least $L/8$. On the right, the two residual constraints for $P_1,P_2$ define a parallelogram of area $4w^2/g$ in the space of witness parameters. At each witness of a retained triple, the allowable completion area is at most $2wL\le 16wg$.}
\label{fig:gaps}
\end{figure}

\paragraph{Counting pairs and triples in sections.}
To estimate how much weight is available, we weight a pair by its horizontal gap and a triple by its horizontal span. For a fixed tube, put
\begin{align*}
  M&=\int_{0<a_1<a_2<1}(a_2-a_1) f(a_1) f(a_2)\diff a_1\diff a_2,\\
  Q&=\int_{0<a_1<a_2<a_3<1}(a_3-a_1)
              f(a_1) f(a_2) f(a_3)\diff a_1\diff a_2\diff a_3.
\end{align*}
To compute the integral of $M$ over all tubes, write $a_i=a(P_i)$ and
use the definition of $f$ to express $M$ as
\[
  M=\frac{1}{(2w)^2}
  \int_{\substack{P_1,P_2\in A\cap T_{u,v}(w)\\ a_1<a_2}}
      (a_2-a_1)\diff P_1\diff P_2.
\]
We can now integrate first over the tubes containing a fixed pair
$P_1,P_2$. The area of their parameter region is $4w^2/(a_2-a_1)$,
so the gap and the normalizing factor $(2w)^{-2}$ give a total
contribution of exactly $1$ per pair. Swapping the two points
interchanges the two possible orders of their first coordinates,
and ties have measure zero. Hence
\[
  \int_{\R^2} M\diff u\diff v=\frac{\alpha^2}{2}.
\]
Also $LM=2Q$. In the product $LM$, the factor $L$ provides a third
horizontal coordinate in addition to the two in $M$. Once the three
coordinates are ordered as $a_1<a_2<a_3$, each of the three choices
of pair in $M$ has the common density $f(a_1)f(a_2)f(a_3)$, multiplied
by the gap of that pair. The three gap factors sum to $2(a_3-a_1)$.
Integrating over the ordered coordinates gives twice the definition
of $Q$. Thus a tube with large $L$ can turn a substantial amount of
pair weight into triple weight.

\paragraph{Keeping sections with enough mass.}
We bound the slopes so that the tubes can later be enlarged to cover whole cells. We also require a lower bound on $L$ so that the retained pair mass yields triple mass via $Q=LM/2$. We restrict the witness parameters to
\[
  W=\{(u,v): |u|\le R,\ L\ge\lambda\},
  \qquad R=\frac{C_0}{\alpha},
  \qquad \lambda=\frac{c_0\alpha}{\log(2/\alpha)},
\]
where $C_0$ is a sufficiently large absolute constant and $c_0$ is a sufficiently small one.

For two points in a tube, subtracting their vertical residuals gives
$|u|(a_2-a_1)\le |b(P_2)-b(P_1)|+2w\le 1+2w$ when $a_1<a_2$.
Thus a tube with $|u|>R$ can contribute only to pairs for which $P_2$
lies within horizontal distance $(1+2w)/R$ to the right of $P_1$.
The area available for $P_2$ is at most $(1+2w)/R$ for each fixed $P_1$,
since the vertical fibers have length at most one. Each such pair
contributes at most $1$ to the integral over these slopes, by the
calculation of the integral of $M$ over all tubes. Integrating over
$P_1\in A$ shows that the discarded pair weight is at most
$\alpha(1+2w)/R$. Choose $C_0$ large enough that this is at most
$\alpha^2/8$.

For a fixed tube, the support of $f$ is contained in an interval of
length at most $h(u)=\min\{1,(1+2w)/|u|\}$, with $h(0)=1$. Indeed,
$f(a)>0$ implies $ua+v\in[-w,1+w]$. In the integral defining $M$,
the gap is therefore at most $h(u)$, and the integral of
$f(a_1)f(a_2)$ over ordered pairs is $L^2/2$. This gives
$M\le h(u)L^2/2$. For each fixed $u$, integrating the definition of
$L$ over $v$ gives $\int_\R L\diff v=\alpha$. For fixed $a$ and
offset $t$, the change of variable $b=ua+v+t$ measures the length
of the vertical fiber of $A$ at $a$. Integrating in $a$ gives
$\alpha$, and averaging over $t$ leaves this value. On tubes with
$L<\lambda$ we use $L^2\le\lambda L$. Thus the pair mass with
$|u|\le R$ and $L<\lambda$ is at most
\[
  \frac{\lambda\alpha}{2}\int_{-R}^R h(u)\diff u
    = O\!\left(\lambda\alpha\log(2/\alpha)\right).
\]
The logarithmic factor comes from integrating $h(u)=O(1/|u|)$
over $1\le |u|\le R=C_0/\alpha$. After fixing $C_0$, choose $c_0$
small enough that this loss is also at most $\alpha^2/8$. At least
half of the original pair weight $\alpha^2/2$ therefore lies in $W$.
Using $Q=LM/2$ and $L\ge\lambda$ on $W$ gives
$\int_W Q\diff u\diff v\ge\lambda\alpha^2/8$.

\paragraph{The cost of separation.}
For $(u,v)\in W$, let $Q^*$ be the contribution to $Q$ from triples
satisfying the separation rule. Fix the outer coordinates $a_1<a_3$.
The middle coordinate can violate the first gap condition only in
$(a_1,a_1+L/8)$ and the second only in $(a_3-L/8,a_3)$. As $f\le1$,
each interval has $f$-mass at most $L/8$. To bound the total loss
from either condition, multiply $L/8$ by the weight
$(a_3-a_1)f(a_1)f(a_3)$ of the outer pair and integrate over
$a_1<a_3$. This costs at most $LM/8=Q/4$ per condition, so
$Q^*\ge Q/2$.

\paragraph{Assigning the weights.}
For three points with distinct first coordinates, order them so that $a_1<a_2<a_3$. The weight is the area of the tube parameter region that retains the triple, multiplied by its horizontal span and divided by $8w^2$.
\[
 K(P_1,P_2,P_3)
   =\frac{a_3-a_1}{8w^2}
     \int_W
     \one\!\left\{\begin{array}{l}
        P_1,P_2,P_3\in T_{u,v}(w),\\
        a_2-a_1\ge L/8,\quad a_3-a_2\ge L/8
     \end{array}\right\}
     \diff u\diff v .
\]
Set $K=0$ when two first coordinates coincide. This weight is
symmetric because the order is determined by the first coordinates.
Positive weight gives a witness with $|u|\le R$ and both gaps at
least $L/8\ge\lambda/8$, as required in part (1). The integral in
$K$ measures a subset of the tubes containing the outer pair,
whose full parameter area is $4w^2/(a_3-a_1)$. Consequently $K$
is half the fraction of that area contributed by the retained
witnesses, which gives $K\le 1/2$.

To compute the total weight on triples in $A$ ordered by first
coordinate, integrate first over the triples for each tube.
Once the first coordinates are fixed, the gap conditions involve
only those coordinates and $L$. Integrating in each height $b(P_i)$
gives the length $2wf(a_i)$ of the corresponding vertical section.
Thus the triple integral is $(2w)^3Q^*$, including the span factor
$a_3-a_1$. Dividing by $8w^2$ gives the total weight
$w\int_W Q^*\diff u\diff v\ge w\lambda\alpha^2/16$, proving
part (2) by the definition of $\lambda$.

For the pair bound, fix $P,Q$ with gap $g>0$. Any tube counted
for a completion $Z$ belongs to $\mathcal T(P,Q)$ and has
$L\le 8g$. For each such tube we can bound the integral over
$Z$ by allowing the entire section $A\cap T_{u,v}(w)$ as the
completion region and using that the span factor is at most one.
Integrating over $Z$ first therefore gives
\[
  \int_A K(P,Q,Z)\diff Z
  \le \frac{1}{8w^2}
      \int_{\mathcal T(P,Q)\cap\{L\le 8g\}} 2wL\diff u\diff v
  \le \frac{1}{8w^2}\frac{4w^2}{g}\,16wg
  =8w .
\]
For $g=0$ the weight vanishes.
\end{proof}

\subsection{Rounding the weights to edges}\label{subsec:rounding}

We now transfer the bound for pairs of points to a bound for pairs of cells.

\begin{lemma}\label{lem:round}
For all sufficiently large $n$ and every collection $J$ of $m\ge V/\sqrt n$ cells, $H[J]$ has a subhypergraph with at least $c\,m^3/(n\log n)$ edges and pair degrees at most $2V/n$, for an absolute constant $c>0$.
\end{lemma}

\begin{proof}
The choice $w=1/(100n)$ leaves room to replace each point by its entire cell. Apply Lemma~\ref{lem:weighted} to the union $A$ of the cells of $J$, with $\alpha=m/V$ and this value of $w$. Since $\alpha\ge n^{-1/2}$, the consecutive gaps of a positive-weight triple are at least $c/(\sqrt n\log n)>\eta$ for large $n$. The three points therefore occupy distinct cells. Also $|u|\le C\sqrt n$ for the witness from the lemma. Moving a point within its cell changes $b-ua-v$ by at most $(1+|u|)\eta$. For large $n$,
\[
  w+(1+C\sqrt n)\eta\le 2w<1/n.
\]
Thus the three entire cells lie in $T_{u,v}(1/n)$ and form an edge of $H[J]$.

\paragraph{Averaging over the cells.}
For an unordered triple $E=\{C_1,C_2,C_3\}$ of distinct cells of $J$, define
\[
  \omega_E
  =V^3\int_{C_1\times C_2\times C_3}K(P,Q,Z)\diff P\diff Q\diff Z .
\]
This is the average of $K$ when one point is chosen uniformly from each cell. The weights lie in $[0,1/2]$, and $\omega_E>0$ implies that $E$ is an edge of $H[J]$.

Since $K$ is symmetric and vanishes on triples with a repeated cell, both summing over unordered cell triples and integrating over triples ordered by first coordinate give one sixth of the corresponding total over labeled triples. Hence
\[
  \mu:=\sum_E\omega_E
  =V^3\int_{\substack{P_1,P_2,P_3\in A\\a(P_1)<a(P_2)<a(P_3)}}
          K(P_1,P_2,P_3)\diff P_1\diff P_2\diff P_3
  \ge c\,\frac{m^3}{n\log n},
\]
where we used $\log(2/\alpha)\le\log n$ for large $n$. For a fixed pair of distinct cells $C_1,C_2$, summing over the possible third cells allows us to apply the pointwise pair bound.
\[
  \sum_{C_3\in J\setminus\{C_1,C_2\}}\omega_{\{C_1,C_2,C_3\}}
  \le V^3\int_{C_1\times C_2}
       \left(\int_A K(P,Q,Z)\diff Z\right)\diff P\diff Q
  \le 8wV .
\]
The last inequality uses $\area(C_1)=\area(C_2)=1/V$.

\paragraph{Selecting edges.}
Keep each edge independently with probability $\omega_E$. The expected number of retained edges is $\mu$, and the expected degree of each pair is at most $8wV$. The Chernoff bound gives probability at most $e^{-\mu/8}$ that fewer than $\mu/2$ edges are retained. It also gives probability at most $\exp(-c'V/n)$ that a fixed pair has degree greater than $2V/n$, for an absolute $c'>0$. Since there are fewer than $V^2$ pairs, the union of these failure events has probability less than one for large $n$, uniformly for $m\ge V/\sqrt n$. Fix a successful subhypergraph for the given $J$.
\end{proof}

\subsection{Controlling vertex degrees as well}\label{subsec:vertex-degrees}

We use the candidate edges from Lemma~\ref{lem:round} to build $H_J$
while keeping the degree of each cell below the bound in (i).

\begin{proof}[Proof of Proposition~\ref{prop:balanced}]
Fix $J$ with $m\ge 2V/\sqrt n$, and put $t=\lfloor c_1m^3/(n\log n)\rfloor$. For every $U\subseteq J$ with $|U|\ge 3m/4$, we have $|U|\ge V/\sqrt n$. Lemma~\ref{lem:round} gives a subhypergraph of $H[U]$ with at least $c(3m/4)^3/(n\log n)$ edges. Choose $c_1>0$ small enough that this is at least $20t$.

We build $H_J$ by adding edges until it has $t$ edges. At each step we obtain a pool of candidates on the cells whose degrees are still below their threshold. A candidate is eligible if it is unused and all three of its pairs are below their threshold.

While fewer than $t$ edges have been chosen, call a vertex saturated if its current degree is at least $12t/m$, and a pair saturated if its current degree is at least $6V/n$. The sum of the vertex degrees is less than $3t$, so fewer than $m/4$ vertices are saturated. Let $U$ consist of the other vertices. By Lemma~\ref{lem:round} we can take a pool of at least $20t$ candidate edges on $U$ in which each pair occurs at most $2V/n$ times.

The sum of the current pair degrees is less than $3t$, so fewer than $3t/(6V/n)$ pairs are saturated. Each rules out at most $2V/n$ candidate edges, and together they rule out at most $t$. Fewer than $t$ candidates have already been used. Thus an unused candidate remains that avoids all saturated pairs, and its vertices all lie in $U$. Add this edge and repeat.

A degree can overshoot its threshold by at most one. When the construction reaches $t$ edges, the maximum vertex degree is at most $12t/m+1$ and the maximum pair degree is at most $6V/n+1$. For large $n$ we have $t/m\ge 1$ and $t\ge c_1m^3/(2n\log n)$, which give the claimed bounds.
\end{proof}

\subsection{Constructing the containers}\label{sec:containers}

We construct $\mathcal C$ by starting with all the cells and applying
the container lemma to $H_J$ at each large container $J$. Each child
retains at most one quarter of the edges of $H_J$. By the single-cell
degree bound, it must therefore omit at least a fixed fraction of the
cells of $J$. This bounds the depth of the tree. Together with the
bound on the number of children at each step, it controls the number
of leaves, which form $\mathcal C$.

We use a fixed edge-reduction form of the container theorem. The maximum degrees refer to hypergraphs without repeated edges.

\begin{lemma}[Saxton--Thomason~\cite{ST15}]\label{lem:container}
There are absolute constants $c_*,C_*,\rho_0>0$ with the following property. Let $G$ be a $3$-uniform hypergraph on $m$ vertices with average degree $d=3e(G)/m>0$. If $0<\rho\le\rho_0$ and
\[
  \frac{\Delta_2(G)}{d\rho}+\frac{1}{d\rho^2}\le c_*,
\]
then there is a family $\mathcal D$ covering the independent sets of $G$ such that every $D\in\mathcal D$ satisfies $e(G[D])\le e(G)/4$ and $\log|\mathcal D|\le C_*m\rho\log(1/\rho)$.
\end{lemma}

This follows from Definition~3.2 and Corollary~3.6 of~\cite{ST15} with edge reduction parameter $1/4$. For uniformity three their co-degree function is $4\delta_2+2\delta_3$, where $\delta_j\le\Delta_j(G)/(d\rho^{j-1})$. Since $\Delta_3(G)\le 1$, the displayed condition makes this co-degree function sufficiently small.

\begin{proof}[Proof of Proposition~\ref{prop:containers}]
Start with the collection of all cells. At a current collection $J$ with $m=|J|$, stop if $m\le A V\sqrt{\log n/n}$, where $A$ is a sufficiently large absolute constant. At a nonterminal collection, $m\ge 2V/\sqrt n$ for large $n$, so Proposition~\ref{prop:balanced} gives $H_J$ with
\[
  d(H_J)\ge c\,\frac{m^2}{n\log n},
  \qquad
  \Delta_2(H_J)\le C\,\frac{V}{n}.
\]
Set $\rho=\kappa\sqrt{n\log n}/m$. This choice makes the second term of the container condition at most $C'/\kappa^2$. The first term is at most
\[
  \frac{C' V\sqrt{\log n}}{\kappa m\sqrt n}
  \le \frac{C'}{\kappa A}.
\]
Choose $\kappa$ and then $A$ sufficiently large. The stopping threshold also gives $\rho\le (\kappa/A)n/V$, which tends to zero. Thus Lemma~\ref{lem:container} applies. Use its containers as the children of $J$ in a tree.

For each child $D$, the lemma gives $e(H_J[D])\le e(H_J)/4$, so at least $3e(H_J)/4$ edges of $H_J$ are lost. Each removed cell meets at most $\Delta_1(H_J)\le C e(H_J)/m$ edges, so $|J\setminus D|\ge 3m/(4C)$. Consequently every nonterminal step reduces the number of cells by a fixed fraction. The tree has depth $O(\log V)=O(\log n)$.

At each nonterminal node, $m\rho=\kappa\sqrt{n\log n}$ and $\log(1/\rho)\le\log V$ for large $n$. The logarithm of the number of children is therefore $O(\sqrt n\,(\log n)^{3/2})$. Multiplying this bound by the depth bounds the logarithm of the number of leaves by $O(\sqrt n\,(\log n)^{5/2})$. Every independent set of $H$ contained in $J$ is also independent in $H_J$, so it is contained in a child of $J$. Following children until the stopping condition holds shows that the leaves cover the independent sets of $H$.

All auxiliary choices are made for subsets of a fixed finite grid. We can fix the entire tree before drawing the coefficient points. Its leaves give $\mathcal C$ with both desired bounds.
\end{proof}

Proposition~\ref{prop:containers} gives the estimate for $\alpha_3(\mathcal S_n)$ by the sampling argument in Section~\ref{sec:cells}. Together with the depth estimate and truncation in Section~\ref{sec:random}, this completes the proof of \added{the planar case of }Theorem~\ref{thm:main}.

\deleted{\medskip

\paragraph{Open problem.}
A natural question is whether $\HD_d(p,d+1)\ge p^{d-o(1)}$ for every
fixed $d\ge 3$. One possible approach is to replace the planar strips
by $n$ random slabs in $\R^d$ of vertical width $2/n$, with defining
coefficients sampled independently and uniformly from $[0,1]^d$.
Extending the container argument would require a higher-dimensional
analogue of Proposition~\ref{prop:balanced} that supplies many edges
while controlling the degrees of individual cells and the codegrees
of sets of two to $d$ cells. The difficulty is that even when the horizontal projections of several coefficient points
are pairwise well separated, they can be nearly affinely dependent,
leaving a large set of tube parameters that can witness the same tuple.}

\begin{additions}
\section{The lower bound on $\HD_d(p,d+1)$}
\label{sec:HDd(p,d+1)}

In this section we prove the lower bound for the $(p,d+1)$ theorem in $\mathbb{R}^d$. A brief overview of the proof was presented in Section~\ref{sec:overview}.
Throughout this section,
$d$ is fixed, and constants may depend on $d$. As in the plane, we
use a geometric estimate to construct the containers, and then bound
the number of sample points in each container.

Put $m=d-1$ and $Q=[0,1]^m$. For $1\le i\le n$, choose
$P_i=(a_i,b_i)$ independently and uniformly from $Q\times[0,1]$ and set
\[
  S_i=\{(x,y)\in\R^m\times\R: |y-a_i\cdot x-b_i|\le 1/n\},
  \qquad \mathcal S_n=\{S_1,\ldots,S_n\}.
\]
The point $(-u,v)$ belongs to the slab associated with $(a,b)$ exactly
when $(a,b)$ belongs to the tube $T_{u,v}(1/n)$, where
\[
  T_{u,v}(w)=\{(a,b)\in\R^d: |b-u\cdot a-v|\le w\}.
\]
For a family $\mathcal F$, let $\alpha_{d+1}(\mathcal F)$ be the largest
size of a subfamily containing no $d+1$ distinct members with a common
point and let $\depth(\mathcal F)$ be the largest number of members of
$\mathcal F$ through any point. With $\ell=\log(en)$ and
$b_d=\max\{3,d-1/d\}$, we prove that,
with probability tending to one as $n\to\infty$,
\begin{equation}\label{eq:ndsample}
  \alpha_{d+1}(\mathcal S_n)\le C_d n^{1/d}\ell^{b_d},
  \qquad
  \depth(\mathcal S_n)\le C_d\frac{\log n}{\log\log n}.
\end{equation}
These are the analogues of (B) and (A), respectively. We also show how to
preserve the intersections used for (B) when we truncate the slabs.

\subsection{The geometric estimate}
The planar estimate controls the weight of triples through a fixed pair.
Here we need a bound for every size of fixed set from two to $d$. If
the $j$ fixed slope vectors span a $(j-1)$-flat, there are $d-j$
directions perpendicular to that flat. Expanding a determinant in
these directions controls the weight of the completions, including
when the remaining slope vectors approach the fixed flat.

Put $w=1/(100n)$ and
\[
 q_j=d+1-j-\frac{j-1}{m},\qquad 2\le j\le d .
\]
The following proposition is the continuous form of the degree bounds.
An integral with $j$ arguments of a kernel fixed is over its other
$d+1-j$ arguments.
\begin{proposition}\label{prop:ndkernel}
If $A\subseteq Q\times[0,1]$ is Borel with volume $\alpha\ge n^{-m/d}$,
there is, for all large $n$, a symmetric measurable function
$K:A^{d+1}\to[0,1]$ satisfying
\begin{equation}
 \int_{A^{d+1}}K\ge c_d w\alpha^{d+1}/\ell^d, \qquad \mbox{and} \label{eq:ndmass}
\end{equation}
\begin{equation}\label{eq:ndrows}
    \sup_{P_1,\ldots,P_j\in A}\int_{A^{d+1-j}}
 K(P_1,\ldots,P_j,Z_1,\ldots,Z_{d+1-j})\,
 \diff Z_1\cdots\diff Z_{d+1-j}
 \le C_d w\alpha^{q_j},
\end{equation}
for all $2\le j\le d$.

\noindent Every tuple with positive weight lies in a tube of halfwidth $w$ with
$\|u\|_\infty\le n^2$.
\end{proposition}

For $k+1$ slope vectors $a_0,\ldots,a_k$, let $D_k$ denote $k!$ times
the $k$-dimensional volume of their simplex, with $D_0=1$. Equivalently,
$D_k$ is the square root of the Gram determinant of $a_i-a_0$,
$1\le i\le k$.
Adjoining another point multiplies this quantity by its distance to the
previous affine span, provided that the previous simplex is nondegenerate.
For $d+1$ projections let $D_*$ be the maximum of $D_m$ over their
$d$-point faces.

\paragraph{A horizontally spread set.}
Fix a sufficiently small $r=r_d>0$, as specified in the normalization
below. Suppose $E\subseteq Q\times[0,1]$ has volume $\gamma>0$ and that
every slab of halfwidth $r$ around an affine hyperplane of $\R^m$ contains,
in its full vertical fibers, at most $(1-\delta)\gamma$ of $E$, where
$0<\delta\le1/2$. For $\pi=(u,v)$ let $\chi_\pi(a,b)$ be the indicator
of $|b-u\cdot a-v|\le w$, and set
\[
 U_\pi(a)=\frac1{2w}\int \one_E(a,b)\chi_\pi(a,b)\diff b,\quad
 L_\pi=\int_Q U_\pi,\quad
 \theta_\pi=\min\{1,\gamma/L_\pi\},
\]
taking $\theta_\pi=1$ if $L_\pi=0$. In particular $U_\pi\le1$ and
$\theta_\pi L_\pi=\min(\gamma,L_\pi)$.

Given $\kappa=c_{d,r}\delta$, call a tuple admissible when \emph{every}
projected $k$-face satisfies
\begin{equation}\label{eq:ndcut}
 D_k\ge\kappa\gamma^{k/m},\qquad 1\le k\le m .
\end{equation}
For a sufficiently large constant $R=R_{d,r}$ write
$\Omega=\{(u,v):\|u\|_\infty\le R\}$. Our kernel on $E$ is
\begin{equation}\label{eq:ndW}
 W(P_1,\ldots,P_{d+1})=
 \one_{\rm admissible}\frac{D_*}{(2w)^d}
 \int_\Omega \theta_\pi\prod_{i=1}^{d+1}\chi_\pi(P_i)\diff\pi.
\end{equation}
For $d$ coefficient points with projected determinant $D_m>0$, the map
from $(u,v)$ to their $d$ graph heights has absolute Jacobian $D_m$.
The unrestricted parameter volume of the tubes containing them is
$(2w)^d/D_m$. Using a face attaining $D_*$ shows $0\le W\le1$.

\begin{lemma}\label{lem:ndspread}
On a horizontally spread set as above, one can choose $\kappa=c_{d,r}\delta$
so that $\int_{E^{d+1}}W\ge c_d w\delta^d\gamma^{d+1}$. With any
$j$ arguments fixed, its completion integral is at most
$C_d w\gamma^{q_j}/\delta$ for $2\le j\le d$.
\end{lemma}
\begin{proof}
\emph{Mass.}
Call an ordered $d$-point simplex broad if each successive projection
after the first is at distance at least $r$ from the affine span of its
predecessors. Each such span is contained in a hyperplane of $\R^m$, so
after any broad prefix the next point may be chosen from a set of measure
at least $\delta\gamma$. Thus broad tuples in $E^d$ have measure at least
$\delta^m\gamma^d$. Their determinant is at least $r^m$. Every one of
their subfaces has determinant at least a positive constant depending only
on $d,r$, since it can be extended to the full simplex by adding altitudes
of length at most $\sqrt m$.

All feasible slopes for a broad tuple lie in $\Omega$ when $R$ is chosen
large enough. Indeed, its prescribed graph heights lie in $[-w,1+w]$,
while the determinant of the height equations is at least $r^m$.
Integrating the $d$ coefficient heights and using the Jacobian just computed
therefore gives
\begin{equation}\label{eq:ndbroad}
 \int_\Omega\int_{Q^d}\one_{\rm broad}D_m
     \prod_{i=1}^d U_\pi(a_i)\diff a_1\cdots\diff a_d\diff\pi
   =\int_{E^d}\one_{\rm broad}\ \ge\ \delta^m\gamma^d.
\end{equation}
Further, $\int_\Omega L_\pi\diff\pi=(2R)^m\gamma$, by integrating the
intercept first. The part of \eqref{eq:ndbroad} on $L_\pi<\lambda\gamma$
is at most $C_d(\lambda\gamma)^m\int_\Omega L_\pi\diff\pi$, since the
inner integral is bounded by $C_d L_\pi^d$. Choose
$\lambda=c'_{d,r}\delta\le1$ small enough to retain at least half the mass.

Fix a remaining plane and a broad simplex. To violate
\eqref{eq:ndcut} with $k$ of its vertices and a new point $z$, the
projection $z$ must be within distance $C_{d,r}\kappa\gamma^{k/m}$ of the
$(k-1)$-flat spanned by those vertices. The intersection with $Q$ of a
radius-$t$ neighborhood of a $(k-1)$-flat has volume at most
$C_dt^{m-k+1}$, by orthogonal projection onto the flat and its complement.
The forbidden volume for this choice of vertices is consequently at most
\[
 C_{d,r}\kappa^{m-k+1}\gamma^{\,k(m-k+1)/m}
 \le C_{d,r}\kappa\gamma,
 \qquad k(m-k+1)-m=(k-1)(m-k)\ge0.
\]
There are only $O_d(1)$ faces. Choose $c_{d,r}$ small enough that the
broad vertices themselves satisfy all cutoffs and the union of the
forbidden sets has volume at most $\lambda\gamma/2$.
Since $U_\pi\le1$ and $L_\pi\ge\lambda\gamma$, the valid extra points
carry $U_\pi$-mass at least $L_\pi/2$. Integrating all $d+1$ heights in
\eqref{eq:ndW} leaves the factor $2w$ and a product of $d+1$ section
densities. Use $D_*\ge D_m$ for the broad face and
$\theta_\pi L_\pi\ge\lambda\gamma$ to obtain
$\int W\ge 2w(\lambda\gamma/2)(\delta^m\gamma^d/2)$.

\emph{Completions.}
Fix $j$ coefficient points, let $k=j-1$ and let $t=m-k=d-j$.
If their $D_k<\kappa\gamma^{k/m}$ the integral is zero. Otherwise center
orthonormal coordinates at the first fixed projection and take the first
$k$ axes in the span of its differences from the other fixed projections.
Let $h_i\in\R^t$ be the normal component of each of the $t+1$ free
projections. All fixed projections have normal component zero.

For each face of $d$ points, write its augmented determinant using
a row of ones, $k$ tangential coordinate rows, and $t$ normal rows.
Expand along the normal rows. A nonzero normal minor uses only free
columns, and every remaining tangential minor is bounded in terms of $d$.
It follows, including for the maximizing face, that
\begin{equation}\label{eq:ndexpansion}
 D_*\ \le C_d\sum_{\substack{I\subseteq[t+1]\\ |I|=t}}|H_I|,
 \qquad H_I=\det(h_i:i\in I).
\end{equation}
For $t=0$ the empty determinant is one. For example, in $d=4$ a fixed
pair leaves a two-dimensional normal space. The three terms on the
right are the areas of the parallelograms spanned by the three possible
pairs of free normal vectors.

For an index set $I$ in this sum,
the projected $m$-dimensional determinant of the $j$ fixed points and the
$t$ free points of $I$ is \emph{exactly} $D_k|H_I|$. Indeed, subtract the
first fixed column and use the block triangular matrix of differences.
Its diagonal blocks are the tangential differences of the fixed points
and the normal components of the free points in $I$.

Bound a term of \eqref{eq:ndexpansion} in the integral of $W$. There is
one free point outside $I$. Integrate it first, dropping admissibility
for the upper bound. On each plane its integral together with the cap is
$2w\theta_\pi L_\pi\le2w\gamma$. The other $d$ points have witness volume
at most $(2w)^d/(D_k|H_I|)$ when $H_I\ne0$. The term with $H_I=0$ contributes
zero before this estimate. Thus the normal minor cancels, and integrating
the $t$ points in $I$ over $E$ gives
\[
 C_d\frac{2w\gamma}{(2w)^d}\int_{E^t}
 |H_I|\,\frac{(2w)^d}{D_k|H_I|}\ \le\
 C_d\,\frac{w\gamma^{t+1}}{D_k}
 \ \le\ C_d\,\frac{w}{\kappa}\gamma^{\,d+1-j-(j-1)/m}.
\]
In the displayed integral the integrand is taken to be zero on $H_I=0$.
Summation over $I$ proves the row bounds.
\end{proof}

\paragraph{Affine normalization.}
We prove Proposition~\ref{prop:ndkernel} by restricting to a core and
rescaling only the horizontal coordinates. Fix $r>0$ so small that
\[
 \zeta=2r(\sqrt m)^{m-1}<1,\qquad
 Q_0=(2\zeta)^{-1}\ge\max\{2,(\sqrt m)^{m-1}\}.
\]
At stage $i=0,1,\ldots$ use
\[
 \delta_i=\tfrac14\bigl(\ell^{-1}+\zeta^{\,i/m}\bigr).
\]
If the current set satisfies the spread assumption with $\delta_i$, stop.
Otherwise retain its part above a horizontal slab of halfwidth $r$ that
contains more than $1-\delta_i$ of the mass. Its horizontal projection
fits into an orthogonal box with $m-1$ side lengths at most $\sqrt m$ and
one side length at most $2r$. Send this box affinely onto $Q$, keeping
the vertical coordinate fixed. The inverse linear part $B_i$ has norm at
most $\sqrt m$ and determinant in absolute value at most $\zeta$.
All box sides can be taken positive. Since $\delta_i\le1/2$, each step
increases the normalized volume by a factor at least $Q_0$. It follows that
there are at most $k_0\le\log(1/\alpha)/\log Q_0=O_d(\ell)$ steps.

Let $\beta$ be the volume retained in the original coordinates and $D$ the
absolute determinant of the composite inverse horizontal map $B$. Then the
terminal normalized set $E$ has volume $\gamma=\beta/D$, and
\[
 c_d\alpha\le\beta\le\alpha,\qquad \beta\le D\le \zeta^{k_0}\le1 .
\]
Indeed, $\sum_{i<k_0}\delta_i=O_d(1)$ by the geometric series and the bound
on $k_0$, hence $\prod_{i<k_0}(1-\delta_i)\ge c_d$. Apply
Lemma~\ref{lem:ndspread} with $\delta=\delta_{k_0}$. Pull its kernel back to
the core and extend it by zero to $A^{d+1}$. Its mass gains a factor
$D^{d+1}$, giving at least $c_dw\delta_{k_0}^d\beta^{d+1}$. As
$\delta_{k_0}\ge1/(4\ell)$ this proves \eqref{eq:ndmass}. A $j$-point
completion integral gains a factor $D^{d+1-j}$ and is at most
\[
 C_d w\beta^{q_j}\,\frac{D^{(j-1)/m}}{\delta_{k_0}}
 \ \le\ 4C_d w\alpha^{q_j}.
\]
Here $q_j\ge0$ and $D^{(j-1)/m}\le \zeta^{k_0/m}$. The decreasing
thresholds use the contraction of the horizontal determinant to absorb the
factor $1/\delta_{k_0}$ simultaneously in every row.

The tube halfwidth is unchanged. A normalized slope $u'$ pulls back to
$B^{-\mathsf T}u'$, and the singular values give
\[
 \|B^{-\mathsf T}\|_2\le\|B\|_2^{m-1}/D
 \le (\sqrt m)^{(m-1)k_0}/\beta\le C_d\alpha^{-2}.
\]
Since $\alpha\ge n^{-m/d}$ and $2m/d<2$, this bound times the fixed
normalized slope bound is at most $n^2$ for large $n$. This completes
the proof of Proposition~\ref{prop:ndkernel}.

\subsection{The cell hypergraph and its containers}

We now pass from weighted tuples to edges on the cells. For each large
set $J$ of cells we select an auxiliary hypergraph $H_J$, whose degree
bounds make it possible to shrink $J$ at each container step.

Write $s_0=d+1$ for the edge size and $L=\ell^d$ for the loss in
\eqref{eq:ndmass}. Partition the coefficient cube into $V=n^{5d}$ cells of
side $\eta=n^{-5}$ using a fixed boundary convention. On these cells define an
$s_0$-uniform hypergraph $H$. Its edges consist of $s_0$ distinct cells
contained together in a single tube of halfwidth $1/n$ with
$\|u\|_\infty\le n^2$. Tuples weighted by Proposition~\ref{prop:ndkernel}
lie in edges when their cells are distinct, since changing all coordinates
of a point within its cell changes the residual by at most
$(1+mn^2)n^{-5}<1/n-w$ for large $n$.

\begin{lemma}\label{lem:ndfinite}
Let $J$ consist of $M=\alpha V$ cells, where $\alpha\ge2n^{-m/d}$. There
is a hypergraph $H_J\subseteq H[J]$ on $J$ with no repeated edges and with
\begin{align*}
 e(H_J)&\ge c_d wM^{s_0}/L,\qquad
 \Delta_1(H_J)\le C_d e(H_J)/M,\\
 \Delta_j(H_J)&\le C_d w V^{s_0-j}\alpha^{q_j},\qquad 2\le j\le d .
\end{align*}
Here $\Delta_j$ is the maximum number of edges through $j$ vertices.
\end{lemma}
\begin{proof}
First take any $J$ with $\alpha\ge n^{-m/d}$. Apply Proposition~\ref{prop:ndkernel}
to its union. Give each unordered tuple of distinct cells weight $V^{s_0}$
times the integral of its kernel over the product in any fixed ordering.
All edge weights are at most one. Tuples with two points in the same cell
have total measure at most $\binom{s_0}{2}\alpha^{s_0-1}/V$, which is
$o(w\alpha^{s_0}/L)$ uniformly in this density range. Dividing the ordered
sum by $s_0!$ gives total weight at least $c'_d w M^{s_0}/L$. Integrating
\eqref{eq:ndrows} over the $j$ fixed cells gives maximum total weight
through them at most $T_j=C_d w V^{s_0-j}\alpha^{q_j}$.

Now take $\alpha\ge2n^{-m/d}$ and set $t=\lfloor c''_d wM^{s_0}/L\rfloor$,
with $c''_d$ so small that any subset of at least $3M/4$ vertices supplies
weights of total at least $10t$. At each step of a greedy selection call
vertices saturated at degree $4s_0t/M$, and $j$-sets saturated at degree
$2d\binom{s_0}{j}T_j$ for $2\le j\le d$. As long as fewer than $t$ edges
are chosen, fewer than $M/4$ vertices are saturated. Apply the weighted
bound afresh to the unsaturated vertices. Their weights still satisfy
the original bounds $T_j$, since their density is at most $\alpha$ and
$q_j\ge0$. There are at most $t/(2dT_j)$ saturated $j$-sets, by summing
their degrees. Avoiding all such sets over all $j$ therefore excludes at
most $t/2$ in weight. Edges already used cost at most another $t$, so a
new edge can be added. After $t$ steps the overshoot in each degree bound is
at most one, absorbed since $t/M,T_j\ge1$ throughout the range for large $n$.
\end{proof}

\begin{proposition}\label{prop:ndcontainers}
The independent sets of $H$ are covered by a deterministic family
$\mathcal C$ of vertex subsets with
\[
 \max_{C\in\mathcal C}|C|/V\le n^{-m/d}\ell^{b_d},\qquad
 \log|\mathcal C|\le C_d n^{1/d}\ell^{b_d}.
\]
\end{proposition}
\begin{proof}
We use the fixed-edge-reduction case of the Saxton--Thomason theorem
\cite[Definition~3.2 and Corollary~3.6]{ST15}. For fixed $s_0$, a hypergraph
of average degree $\bar d$ on $M$ vertices has containers spanning at most
one quarter of its edges when $\rho>0$ and
all the ratios $\Delta_j/(\bar d\rho^{j-1})$, $2\le j\le s_0$, and
the parameter $\rho$ are sufficiently small in terms of $s_0$. Their number has logarithm at most
$C_{s_0}M\rho\log(1/\rho)$.

Starting at all cells, stop a node of density $\alpha$ once
$\alpha\le\alpha_0=n^{-m/d}\ell^{b_d}$. Otherwise take $H_J$ from
Lemma~\ref{lem:ndfinite}, and choose
\[
 s(\alpha)=K_d\bigl((nL)^{1/d}+L\alpha^{-1/m}\bigr),\qquad
 \rho=s(\alpha)/M.
\]
Here $K_d$ is a sufficiently large constant. We have $\rho=o(1)$
uniformly at these nodes, since $V=n^{5d}$ and $\alpha>\alpha_0$.
The mean degree is at least $c_d w M^d/L$. For $2\le j\le d$ the ratios
in the container condition are at most
\[
 \frac{C_d L}{(s(\alpha)\alpha^{1/m})^{j-1}}
 \le \frac{C_d}{K_d^{j-1}L^{j-2}} .
\]
For $j=s_0$ the ratio is at most $C_dnL/s(\alpha)^d$, since
$\Delta_{s_0}=1$. This verifies all conditions.

Use the containers of $H_J$ as children of the current node. Each misses at
least $3e(H_J)/4$ edges. The vertex degree bound thus makes each child at most
$\lambda_d M$ in size for a constant $\lambda_d<1$. All independent sets of
$H[J]$ are covered by the children because $H_J\subseteq H[J]$. The depth
is $O_d(\ell)$. To count the leaves when $s(\alpha)$ varies, consider any
path whose internal densities in order are $\alpha_1,\ldots,\alpha_h$.
They decrease by a factor at most $\lambda_d$ per step, so
\[
 \sum_{i=1}^h\alpha_i^{-1/m}
 \le\alpha_h^{-1/m}\sum_{t\ge0}\lambda_d^{t/m}
 \le C_d\alpha_0^{-1/m}.
\]
Since $\log(1/\rho)\le\log V$, the sum of the logarithms of the branching
numbers on this path is at most
\[
 C_d\bigl[\ell^2(nL)^{1/d}+\ell L\alpha_0^{-1/m}\bigr]
 \le C'_d n^{1/d}\ell^{b_d}.
\]
For the last bound the two conditions on $b_d$ are $b_d\ge3$ and
$b_d\ge d+1-b_d/m$, the latter equivalent to $b_d\ge d-1/d$.
Assign equal probabilities to all children at each branch. The bound on
every path gives each leaf probability at least
$\exp(-C'_dn^{1/d}\ell^{b_d})$, proving the asserted count. All choices
are on a finite grid and can be fixed before the sample is drawn.
\end{proof}

\subsection{Sampling and the piercing bound}

The containers give the first bound in \eqref{eq:ndsample}. For a fixed terminal container its occupancy by sampled labels is binomial
with mean at most $n\alpha_0=n^{1/d}\ell^{b_d}$. Choose
$h=\lceil B_dn^{1/d}\ell^{b_d}\rceil$. The binomial upper tail and the
union bound give
\[
 \Prob(\text{some terminal container receives at least $h$ labels})
 \le \exp(C_d n^{1/d}\ell^{b_d})(e/B_d)^h=o(1)
\]
for sufficiently large $B_d$. If a subfamily has no intersecting
$(d+1)$-tuple, its occupied cells are independent in $H$: representatives
of any edge would all have the common point $(-u,v)$. This observation also
applies when sample cells repeat. Such a subfamily is therefore contained
in one terminal container and has fewer than $h$ labels with high probability.

\paragraph{Depth.}
At a fixed point the probability of lying in one random slab
is at most $2/n$, by conditioning on $a_i$ and integrating $b_i$.
Almost surely any $d$ of the slab normals $(-a_i,1)$ are linearly independent.
Their slabs have bounded intersection. A nonempty intersection of at
least $d$ slabs thus has a vertex supported by $d$ independent signed
boundary planes from distinct slabs. Conditional on the $d$ defining
slabs, each remaining slab contains this vertex independently with
probability at most $2/n$. A union bound over the candidates gives
\[
 \Prob\bigl(\depth(\mathcal S_n)\ge t\bigr)
 \le 2^d\binom nd\left(\frac{2e}{t-d}\right)^{t-d}
 \quad\text{for integers }t>d .
\]
Taking $t=r_n=\lceil(2d+2)\log n/\log\log n\rceil$ proves the second
part of \eqref{eq:ndsample}.

\paragraph{Compactness and the choice of $n$.}
Intersect each slab $S_i$ with $[-n^2,n^2]^m\times\R$ and call the
resulting body $F_i$. It is the affine image of $[-n^2,n^2]^m\times[-1/n,1/n]$ by
$(x,z)\mapsto(x,a_i\cdot x+b_i+z)$, hence is a full-dimensional compact
parallelotope. The bodies are distinct almost surely, for example by their
vertical sections at $x=0$. Every cell edge still certifies a common point
of the truncated slabs, since its slope satisfies $\|u\|_\infty\le n^2$.
Thus the same container occupancy bound proves the first part of
\eqref{eq:ndsample} for the compact family. Truncation can only decrease depth.

For a given sufficiently large integer $p$, take
$n=\lfloor c'_dp^d/(\log p)^{db_d}\rfloor$ with a sufficiently small
$c'_d>0$. Then $h\le p$ and $n\ge p$. Fix a sample satisfying both
estimates and let $\mathcal F_n=\{F_1,\ldots,F_n\}$. It has the
$(p,d+1)$ property and
\[
  \tau(\mathcal F_n)\ge \frac{n}{r_n}
     =\Omega_d\!\left(\frac{n\log\log n}{\log n}\right).
\]
Since $db_d+1=\max\{3d+1,d^2\}$ and $\log n\asymp_d\log p$, this gives
Theorem~\ref{thm:main} for $d\ge 3$.
\end{additions}

\section{The lower bound on \texorpdfstring{$\HD_d(p,q)$}{HDd(p,q)}}
\label{sec:general}

We prove the second bound in Theorem~\ref{thm:main} for each fixed
$q\ge d+2$. We use the random slabs $\mathcal S_n$ of
Section~\ref{sec:HDd(p,d+1)}, whose depth is
$O_d(\log n/\log\log n)$ with probability tending to one by
\eqref{eq:ndsample}. Constants may depend on $d$ and $q$. As in that
section, put $m=d-1$, $Q=[0,1]^m$, $\ell=\log(en)$ and $w=1/(100n)$,
and write $s=q-d$. To establish the $(p,q)$ property, we need a
weighted count of $q$-tuples lying in a tube $T_{u,v}(w)$.
The factor $w$ in \eqref{eq:ndmass} becomes $w^s$, since $d$
coefficient points specify a candidate hyperplane and $s$ further
points must lie in its tube.

\subsection{The geometric estimate}\label{subsec:fixed-q-geometry}

For $2\le j<q$, use the exponents
\[
 (\nu_j,q_j)=
 \begin{cases}
  (s,\ q-j-(j-1)/m),& 2\le j\le d,\\
  (q-j,\ q-j-1),& d<j<q .
 \end{cases}
\]
Both exponents are nonnegative.

\begin{proposition}\label{prop:fixed-q-kernel}
Let $A\subseteq Q\times[0,1]$ be a Borel set of volume
$\alpha\ge n^{-m/(q-1)}$. For all sufficiently large $n$ there is a
symmetric measurable function $K:A^q\to[0,1]$, supported on tuples in
tubes $T_{u,v}(w)$ with $\|u\|_\infty\le n^2$, such that
\begin{align}
 \int_{A^q}K &\ge c\,w^s\alpha^q/\ell^{q-1},
                 \label{eq:fixed-q-mass}\\
 \sup_{P_1,\ldots,P_j\in A}\int_{A^{q-j}}
   K(P_1,\ldots,P_j,Z_1,\ldots,Z_{q-j})\,\diff Z_1\cdots\diff Z_{q-j}
 &\le Cw^{\nu_j}\alpha^{q_j},\qquad 2\le j<q .
                 \label{eq:fixed-q-rows}
\end{align}
\end{proposition}

We use a determinant estimate to select the $s$ additional points
simultaneously.

\begin{lemma}\label{lem:fixed-q-determinants}
Let $Y_1,\ldots,Y_t$ be independent random vectors in $\R^r$, with
$t\le r$ and densities bounded by $B$. For $z>0$,
\[
 \Prob\bigl(\lVert Y_1\wedge\cdots\wedge Y_t\rVert\le z\bigr)
 \le C_{r,t} B z^{r/t}.
\]
\end{lemma}

\begin{proof}
For $t=1$ this is the volume bound for a ball. For $t>1$, order the
vectors by choosing at each step a vector at greatest distance from
the span of those already chosen. Consider one possible order and
condition on its first $t-1$ vectors. On the event that this order is
chosen, their successive distances are $b_1\ge\cdots\ge b_{t-1}>0$
outside a null set. In the corresponding Gram--Schmidt directions,
the last vector has coordinates of absolute value at most $b_i$.
Its distance to their span is at most
\[
 u=\min\{b_{t-1},\,z/(b_1\cdots b_{t-1})\}.
\]
The allowed volume is at most $C_{r,t}(b_1\cdots b_{t-1})u^{r-t+1}$.
Since $u^t\le(b_1\cdots b_{t-1})u\le z$, this is at most
$C_{r,t}z^{r/t}$. Independence gives the asserted bound after
summing over the possible orders.
\end{proof}

\begin{proof}[Proof of Proposition~\ref{prop:fixed-q-kernel}]
\emph{A horizontally spread set.}
Work first on a set $E$ of volume $\gamma$ with spread parameter
$\delta$ as in Section~\ref{sec:HDd(p,d+1)}. Use the same
$U_\pi,L_\pi,\theta_\pi,\Omega$ as in \eqref{eq:ndW}. For a $q$-tuple,
let $D_*$ be the largest $D_m$ of its $d$-point faces. Require
\eqref{eq:ndcut} on every projected face with
$\kappa=c_{d,q,r}\delta$, and set
\[
 W(P_1,\ldots,P_q)=
 \one_{\rm admissible}\frac{D_*}{(2w)^d}
 \int_\Omega\theta_\pi^s\prod_{i=1}^q\chi_\pi(P_i)\diff\pi .
\]
The height Jacobian for a face attaining $D_*$ again gives $0\le W\le1$.

\emph{Mass.}
Begin with a broad $d$-tuple as in the proof of
Lemma~\ref{lem:ndspread}. By the argument for \eqref{eq:ndbroad},
restricting to $L_\pi\ge c_d\delta\gamma$ retains at least half of the
broad-tuple contribution. For each such tube and broad tuple, draw the
$s$ new projections independently with density $U_\pi/L_\pi$.

Consider a projected $k$-face having $b\ge1$ broad vertices and
$t=k+1-b\ge1$ new vertices. Project the new differences from a broad
vertex onto the orthogonal complement of the span of the broad
differences. Its dimension is $r'=m-b+1$. The projected densities are
bounded by $C_d/L_\pi$, since the original densities are bounded by
$1/L_\pi$ and supported in $Q$. The face determinant is the product of
the broad base determinant, bounded below by a positive constant,
and the exterior product of the $t$ projected differences. By
Lemma~\ref{lem:fixed-q-determinants}, the probability that
\eqref{eq:ndcut} fails on this face is at most
\[
 C\frac{\kappa^{r'/t}\gamma^{kr'/(mt)}}{L_\pi}
 \le C\frac{\kappa}{\delta},
\]
where $r'\ge t$ and $kr'-mt=(b-1)(m-k)\ge0$. For a face with only
new vertices, condition on its first vertex and use the same bound
with $b=1$ and $t=k$. Faces with only broad vertices already satisfy
the cutoffs. Choosing $c_{d,q,r}$ sufficiently small makes the sum
of all failure probabilities at most $1/2$.

After integrating the heights, the $s$ new points thus contribute
at least a constant times $(w\theta_\pi L_\pi)^s$. On the retained
tubes $\theta_\pi L_\pi=\min(\gamma,L_\pi)\ge c_d\delta\gamma$.
Multiplying this by the broad contribution from \eqref{eq:ndbroad}
and using $D_*\ge D_m$ for the broad face gives
$\int_{E^q}W\ge c w^s\delta^{q-1}\gamma^q$.

\emph{Completions.}
For $2\le j\le d$, the fixed points have projected determinant
$D_{j-1}\ge\kappa\gamma^{(j-1)/m}$ whenever the weight is positive.
Expand the $d$-point determinants defining $D_*$ in the $d-j$
normal directions to the fixed face, as in the proof of
Lemma~\ref{lem:ndspread}. This bounds $D_*$ by a constant times the
sum of the absolute normal minors, one for each choice of $d-j$
free points. In the term for a minor $H_I$, the corresponding
$d$-point determinant is $D_{j-1}|H_I|$. Integrating the other $s$
points first contributes at most $(2w\gamma)^s$, by the cap
$\theta_\pi^s$. The tube parameter volume for the $d$ points is
at most $(2w)^d/(D_{j-1}|H_I|)$ when $H_I\ne0$. The minor cancels
and integrating the $d-j$ chosen free points gives a bound of
$Cw^s\gamma^{q-j}/D_{j-1}$. Thus the row integral is at most
$Cw^s\gamma^{q_j}/\delta$.

For $d<j<q$, write $e=q-j$. Integrating all the free points first
costs at most $(2w\gamma)^e$, since $e<s$ and
$\theta_\pi^s L_\pi^e\le\gamma^e$. Any $d$ of the fixed points have
determinant at least $\kappa\gamma$ on admissible tuples, so they
bound the tube parameter volume by $(2w)^d/(\kappa\gamma)$.
Using $D_*\le C_d$ gives $Cw^e\gamma^{e-1}/\delta$.

\emph{Affine normalization.}
Use the horizontal rescaling from the proof of
Proposition~\ref{prop:ndkernel}. In its notation, the retained volume
is $\beta\asymp_d\alpha$, the determinant of the inverse horizontal
map satisfies $D\le\zeta^{k_0}$, and the final spread parameter is
$\delta_{k_0}=(\ell^{-1}+\zeta^{k_0/m})/4$. The normalized volume
is $\gamma=\beta/D$. Pulling back multiplies the mass by $D^q$,
and gives \eqref{eq:fixed-q-mass}. A $j$-point row gains a factor
$D^{q-j}$. In both ranges $q-j-q_j\ge1/m$, so
$D^{q-j-q_j}/\delta_{k_0}\le4$. Since $q_j\ge0$ and $\beta\le\alpha$,
this proves \eqref{eq:fixed-q-rows}. The slope bound follows from
$\lVert B^{-\mathsf T}\rVert_2\le C_d\alpha^{-2}$ for the inverse
horizontal map $B$, as in Section~\ref{sec:HDd(p,d+1)}, using
$2m/(q-1)<2$.
\end{proof}

\subsection{Containers and sampling}\label{subsec:fixed-q-containers}

On the $V=n^{5dq}$ cells of side $n^{-5q}$, define $H$ by taking
as edges all $q$-sets of cells in a common tube of halfwidth $1/n$
with $\|u\|_\infty\le n^2$. A positively weighted tuple on distinct cells
gives an edge, since the residual variation in a cell is at most
$(1+mn^2)n^{-5q}<1/n-w$.

Write $L=\ell^{q-1}$. Let $J$ consist of $M=\alpha V$ cells, with
$\alpha\ge2n^{-m/(q-1)}$. The rounding argument of
Lemma~\ref{lem:ndfinite}, applied to Proposition~\ref{prop:fixed-q-kernel},
gives a simple $H_J\subseteq H[J]$ such that
\begin{align*}
 e(H_J)&\ge c w^s M^q/L,\qquad
 \Delta_1(H_J)\le C e(H_J)/M,\\
 \Delta_j(H_J)&\le Cw^{\nu_j}V^{q-j}\alpha^{q_j},
                           \qquad 2\le j<q .
\end{align*}
Indeed, tuples with a repeated cell cost at most
$\binom q2\alpha^{q-1}/V=o(w^s\alpha^q/L)$, even at the smallest
density in Proposition~\ref{prop:fixed-q-kernel}. Integrating over
the product of distinct cells and multiplying by $V^q$ gives
weights at most one and a row bound
$T_j=Cw^{\nu_j}V^{q-j}\alpha^{q_j}$. In the greedy selection of
$t=\lfloor c'w^s M^q/L\rfloor$ edges, use saturation thresholds
$4qt/M$ on vertices and $2q\binom qjT_j$ on $j$-sets. At each
step at least $3M/4$ vertices remain unsaturated. Recompute the
weights on those vertices, whose density is still at least
$n^{-m/(q-1)}$. Their total weight is at least $10t$ for small
enough $c'$, and their row bounds remain $T_j$ since $q_j\ge0$.
Summing the $j$-degrees shows that there are at most $t/(2qT_j)$
saturated $j$-sets. They cost at most $t/2$ in weight summed over
all $j$, and already chosen edges cost at most $t$.
Thus the selection continues to $t$ edges. The integer overshoots
are absorbed in the constants since $t/M,T_j$ tend to infinity
on this grid.

Set $\xi=s/(q-1)$. Stop the container recursion when the density
is at most $\alpha_0=n^{-m/(q-1)}\ell^3$. At an internal node use
$H_J$ and the Saxton--Thomason parameter $\rho=S/M$, where
$S=K_0 n^\xi\ell$ for a sufficiently large constant $K_0$.
The mean degree $\bar d$ of $H_J$ is at least $c w^sM^{q-1}/L$.
Substituting the codegree bounds shows
\[
 \frac{\Delta_j(H_J)}{\bar d\,\rho^{j-1}}
 \le
 \begin{cases}
  C L/(S\alpha^{1/m})^{j-1}, & 2\le j\le d,\\
  C L n^{j-d}/(\alpha S^{j-1}), & d<j<q,\\
  C L n^s/S^{q-1}, & j=q .
 \end{cases}
\]
In the first range $S\alpha^{1/m}$ grows at least as a positive
power of $n$, because $s\ge2$. In the second range substitution
of $\alpha_0$ bounds the row by
$C K_0^{1-j}n^{-m(q-1-j)/(q-1)}\ell^{q-j-3}$, which tends to zero
also when $j=q-1$. The last ratio is at most $C/K_0^{q-1}$, and
$\rho=o(1)$ on the chosen grid.

Apply Saxton--Thomason~\cite[Definition~3.2 and Corollary~3.6]{ST15}
as in Section~\ref{sec:HDd(p,d+1)}. Each child omits at least three
quarters of the edges of $H_J$, so the bound on $\Delta_1$ gives
a fixed-factor decrease in size. The independent sets of $H$ are
covered by the terminal containers $\mathcal C$. There are $O(\ell)$
levels and each node has $\exp(O(S\log V))$ children. Thus
$\log|\mathcal C|=O(n^\xi\ell^3)$ and each terminal container has
density at most $\alpha_0$.

For a fixed terminal container the occupancy by sampled labels is
binomial with mean at most $n\alpha_0=n^\xi\ell^3$. The binomial tail
and a union bound show that, with probability tending to one, every
terminal container receives fewer than
$p_n=\lceil B_0 n^\xi\ell^3\rceil$ labels for a sufficiently large
constant $B_0$. This counts labels even when cells repeat.

Truncate at $F_i=S_i\cap([-n^2,n^2]^m\times\R)$ as in
Section~\ref{sec:HDd(p,d+1)}. Edge witnesses survive, so any subfamily
without an intersecting $q$-tuple occupies independent cells and
has fewer than $p_n$ members. The distinct bodies $F_i$ require at
least $c_dn\log\log n/\log n$ piercing points by the depth bound.
For each sufficiently large $p$, take
\[
 n=\left\lfloor c'_{d,q}\,
       \frac{p^{1/\xi}}{(\log p)^{3d}}\right\rfloor
\]
with $c'_{d,q}>0$ sufficiently small. Since $1/\xi=(q-1)/(q-d)\le d$
and $\ell\asymp_{d,q}\log p$, the power $3d$ in the denominator
ensures that $p_n\le p$. Also $1/\xi>1$, so $p\le n$ for all
sufficiently large $p$. The truncated family therefore satisfies
the $(p,q)$ property, and substituting this choice of $n$ into
the piercing bound yields
\[
 \HD_d(p,q)\ge c_{d,q}\,
 \frac{p^{(q-1)/(q-d)}\log\log p}{(\log p)^{3d+1}}.
\]
This proves the second bound in Theorem~\ref{thm:main}.

\medskip

\noindent\textbf{AI disclosure.} The construction was developed with the assistance of ChatGPT 6.2 Astra, ChatGPT 6.2 Sol and ChatGPT 6.1 Sol. The author rewrote the main part of the paper and edited the rest of the proof.

\end{document}